# Numerical Solutions of 2-D Steady Incompressible Driven Cavity Flow at High Reynolds Numbers

*E. Erturk[1] and T.C. Corke[2]
and
C. Gökçöl[1]

[1] *Energy Systems Engineering Department, Gebze Institute of Technology,
Gebze, Kocaeli 41400, Turkey*
[2] *Department of Aerospace and Mechanical Engineering,
University of Notre Dame, IN 46556, USA*

SUMMARY

Numerical calculations of the 2-D steady incompressible driven cavity flow are presented. The Navier-Stokes equations in streamfunction and vorticity formulation are solved numerically using a fine uniform grid mesh of $601 \times 601$. The steady driven cavity solutions are computed for Re $\leq$ 21,000 with a maximum absolute residuals of the governing equations that were less than $10^{-10}$. A new quaternary vortex at the bottom left corner and a new tertiary vortex at the top left corner of the cavity are observed in the flow field as the Reynolds number increases. Detailed results are presented and comparisons are made with benchmark solutions found in the literature.

## 1. INTRODUCTION

Numerical methods for 2-D steady incompressible Navier-Stokes (N-S) equations are often tested for code validation, on a very well known benchmark problem; the lid-driven cavity flow. Due to the simplicity of the cavity geometry, applying a numerical method on this flow problem in terms of coding is quite easy and straight forward. Despite its simple geometry, the driven cavity flow retains a rich fluid flow physics manifested by multiple counter rotating recirculating regions on the corners of the cavity depending on the Reynolds number. In the literature, it is possible to find different numerical approaches which have been applied to the driven cavity

*Correspondence to: ercanerturk@gyte.edu.tr

Download figures, tables, data files, fortran codes and etc. from `http://www.cavityflow.com`. All the numerical solutions presented in this study is available to public in this web site.

Contract/grant sponsor: GYTE-BAP; contract/grant number: BAP-2003-A-22

flow problem [1-34]. Though this flow problem has been numerically studied extensively, still there are some points which are not agreed upon. For example; 1) an interesting point among many studies is that different numerical solutions of cavity flow yield about the same results for Re $\leq$ 1,000 however, start to deviate from each other for larger Re, 2) also another interesting point is that while some studies predict a periodic flow at a high Reynolds number, some others present steady solutions for even a higher Reynolds number. The objective of this work is then to investigate these two points and present accurate very fine grid numerical solutions of steady 2-D driven cavity flow.

Among the numerous studies found in the literature, we will give a brief survey about some of the significant studies that uses different types of numerical methods on the driven cavity flow. In doing so the emphasis will be given on three points; on the numerical method used, on the spatial order of the numerical solution and on the largest Reynolds number achieved.

Recently, Barragy & Carey [2] have used a $p$-type finite element scheme on a $257 \times 257$ strongly graded and refined element mesh. They have obtained a highly accurate ($\Delta h^8$ order) solutions for steady cavity flow solutions up to Reynolds numbers of Re=12,500. Although Barragy & Carey [2] have presented qualitative solution for Re=16,000 they concluded that their solution for Re=16,000 was under-resolved and needed a greater mesh size.

Botella & Peyret [6] have used a Chebyshev collocation method for the solution of the lid-driven cavity flow. They have used a subtraction method of the leading terms of the asymptotic expansion of the solution of the N-S equations in the vicinity of the corners, where velocity is discontinuous, and obtained a highly accurate spectral solutions for the cavity flow with a maximum of grid mesh of $N = 160$ (polynomial degree) for Reynolds numbers Re $\leq$ 9,000. They stated that their numerical solutions exhibit a periodic behavior beyond this Re.

Schreiber & Keller [28] have introduced an efficient numerical technique for steady viscous incompressible flows. The non-linear differential equations are solved by a sequence of Newton and chord iterations. The linear systems associated with the Newton iteration is solved by $LU$-factorization with partial pivoting. Applying repeated Richardson extrapolation using the solutions obtained on different grid mesh sizes (maximum being $180 \times 180$), they have presented high-order accurate ($\Delta h^8$ order in theory) solutions for Reynolds numbers Re $\leq$ 10,000.

Benjamin & Denny [5] have used a method of relaxing the algebraic equations by means of the ADI method, with a non-uniform iteration parameter. They have solved the cavity flow for Re $\leq$ 10,000 with three different grid mesh sizes (maximum being $101 \times 101$) and used a $\Delta h^n$ extrapolation to obtain the values when $\Delta h \to 0$.

Wright & Gaskell [36] have applied the Block Implicit Multigrid Method (BIMM) to the SMART and QUICK discretizations. They have presented cavity flow results obtained on a $1024 \times 1024$ grid mesh for Re $\leq$ 1,000.

Nishida & Satofuka [24] have presented a new higher order method for simulation of the driven cavity flow. They have discretized the spatial derivatives of the N-S equations using a modified differential quadrature (MDQ) method. They have integrated the resulting system of ODEs in time with fourth order Runge-Kutta-Gill (RKG) scheme. With this they have presented spatially $\Delta h^{10}$ order accurate solutions with grid size of $129 \times 129$ for Re $\leq$ 3,200.

Liao [20], and Liao & Zhu [21], have used a higher order streamfunction-vorticity boundary element method (BEM) formulation for the solution of N-S equations. With this they have presented solutions up to Re=10,000 with a grid mesh of $257 \times 257$.

Hou et. al. [17] have used Lattice Boltzmann Method for simulation of the cavity flow. They have used $256 \times 256$ grid points and presented solutions up to Reynolds numbers of Re=7,500.

Goyon [13] have solved the streamfunction and vorticity equations using Incremental Unknowns. They have presented steady solutions for Re $\leq$ 7,500 on a maximum grid size of $256 \times 256$.

Rubin & Khosla [27] have used the strongly implicit numerical method with $2 \times 2$ coupled streamfunction-vorticity form of the N-S equations. They have obtained solutions for Reynolds numbers in the range of Re $\leq$ 3,000 with a grid mesh of $17 \times 17$.

Ghia et. al. [12] later have applied a multi-grid strategy to the coupled strongly implicit method developed by Rubin & Khosla [27]. They have presented solutions for Reynolds numbers as high as Re=10,000 with meshes consisting of as many as $257 \times 257$ grid points.

Gupta [15] has used a fourth ($\Delta h^4$) order compact scheme for the numerical solution of the driven cavity flow. He has used a 9 point ($3\times3$) stencil in which the streamfunction and vorticity equations are approximated to fourth order accuracy. He has used point-SOR type of iteration and presented steady cavity flow solutions for Re $\leq$ 2,000 with a maximum of $41 \times 41$ grid mesh.

Li et. al. [19] have used a fourth ($\Delta h^4$) order compact scheme which had a faster convergence than that of Gupta [15]. They have solved the cavity flow with a grid size of $129 \times 129$ for Re $\leq$ 7,500.

Bruneau & Jouron, [7] have solved the N-S equations in primitive variables using a full multigrid-full approximation storage (FMG-FAS) method. With a grid size $256 \times 256$, they have obtained steady solutions for Re $\leq$ 5,000.

Grigoriev & Dargush [14] have presented a BEM solution with improved penalty function technique using hexagonal subregions and they have discretized the integral equation for each subregion as in FEM. They have used a non-uniform mesh of 5040 quadrilateral cells. With this they were able to solve driven cavity flow up to Re=5,000.

Aydin & Fenner [1] have used BEM formulation with using central and upwind finite difference scheme for the convective terms. They have stated that their formulation lost its reliability for Reynolds numbers greater than 1,000.

To the authors best knowledge, in all these studies among with other numerous papers found in the literature, the maximum Reynolds number achieved for the *2-D steady incompressible* flow in a lid-driven cavity is Re=12,500 and is reported by Barragy & Carey [2]. Although Barragy & Carey [2] have presented solutions for Re=16,000, their solution for this Reynolds number display oscillations related to boundary layer resolution issues due to a coarse mesh. Nallasamy & Prasad [22] have presented steady solutions for Reynolds numbers up to Re $\leq$ 50,000, however, their solutions are believed to be inaccurate as a result of excessive numerical dissipation caused by their first order upwind difference scheme.

A very brief discussion on computational as well as experimental studies on the lid-driven cavity flow can be found in Shankar & Deshpande [29].

Many factors affect the accuracy of a numerical solution, such as, the number of grids in the computational mesh ($\Delta h$), and the spatial discretization order of the finite difference equations, and also the boundary conditions used in the solution.

It is an obvious statement that as the number of grids in a mesh is increased (smaller $\Delta h$) a numerical solution gets more accurate. In this study the effect of number of grid points in a mesh on the accuracy of the numerical solution of driven cavity flow is investigated, especially as the Reynolds number increases. For this, the governing N-S equations are solved on progressively increasing number of grid points (from $128\times128$ to $600\times600$) and the solutions are compared with the highly accurate benchmark solutions found in the literature.

Order of spatial discretization is another factor on the accuracy of a numerical solution. Finite element and spectral methods provide high spatial accuracy. This study also investigates the effect of spacial accuracy on the numerical solution of driven cavity flow. For this, following [28], repeated Richardson extrapolation is used and then the extrapolated results are compared with the benchmark solutions.

Another factor on the accuracy of a numerical solution is the type and the order of the numerical boundary conditions used in computation. For driven cavity flow, this subject was briefly discussed by Weinan & Jian-Guo [34], Spotz [30], Napolitano et. al. [23] and also by Gupta & Manohar [16].

Until Barragy & Carey [2], many studies have presented steady solutions of driven cavity flow for Re=10,000 ([28], [5], [21] and [12]). They ([2]) have presented solutions for Re=12,500. In all these studies that have presented high Reynolds number numerical solutions, the largest grid size used was $256 \times 256$. In this study, the grid size will be increased up to $600 \times 600$ and the effect of the grid size on the largest computable Reynolds number solution of the driven cavity flow will be investigated.

This paper presents accurate, very fine grid ($600 \times 600$) numerical solutions of 2-D steady incompressible flow in a lid-driven cavity for Reynolds numbers up to $Re = 21,000$. A detailed comparisons of our results with mainly the studies mentioned above as well as with other studies not mentioned here, will be done.

## 2. NUMERICAL METHOD

For two-dimensional and axi-symmetric flows it is convenient to use the streamfunction ($\psi$) and vorticity ($\omega$) formulation of the Navier-Stokes equations. In non-dimensional form, they are given as

$$\frac{\partial^2 \psi}{\partial x^2} + \frac{\partial^2 \psi}{\partial y^2} = -\omega \tag{1}$$

$$\frac{1}{\text{Re}} \left( \frac{\partial^2 \omega}{\partial x^2} + \frac{\partial^2 \omega}{\partial y^2} \right) = \frac{\partial \psi}{\partial y} \frac{\partial \omega}{\partial x} - \frac{\partial \psi}{\partial x} \frac{\partial \omega}{\partial y} \tag{2}$$

where, Re is the Reynolds number, and $x$ and $y$ are the Cartesian coordinates. The N-S equations are non-linear, and therefore need to be solved in an iterative manner. In order to have an iterative numerical algorithm, pseudo time derivatives are assigned to these equations (1 and 2) such as

$$\frac{\partial \psi}{\partial t} = \frac{\partial^2 \psi}{\partial x^2} + \frac{\partial^2 \psi}{\partial y^2} + \omega \tag{3}$$

$$\frac{\partial \omega}{\partial t} = \frac{1}{\text{Re}} \frac{\partial^2 \omega}{\partial x^2} + \frac{1}{\text{Re}} \frac{\partial^2 \omega}{\partial y^2} + \frac{\partial \psi}{\partial x} \frac{\partial \omega}{\partial y} - \frac{\partial \psi}{\partial y} \frac{\partial \omega}{\partial x} \tag{4}$$

Since we are seeking a steady solution, the order of the pseudo time derivatives will not affect the final solution. Therefore, an implicit Euler time step is used for these pseudo time derivatives, which is first order ($\Delta t$) accurate. Also for the non-linear terms in the vorticity equation, a first order ($\Delta t$) accurate approximation is used and the governing equations (3 and 4) are written as the following

$$\psi^{n+1} - \Delta t \frac{\partial^2 \psi}{\partial x^2}^{n+1} - \Delta t \frac{\partial^2 \psi}{\partial y^2}^{n+1} = \psi^n + \Delta t \omega^n \tag{5}$$

$$\omega^{n+1} - \Delta t \frac{1}{\text{Re}} \frac{\partial^2 \omega}{\partial x^2}^{n+1} - \Delta t \frac{1}{\text{Re}} \frac{\partial^2 \omega}{\partial y^2}^{n+1} - \Delta t \frac{\partial \psi}{\partial x}^n \frac{\partial \omega}{\partial y}^{n+1} + \Delta t \frac{\partial \psi}{\partial y}^n \frac{\partial \omega}{\partial x}^{n+1} = \omega^n \tag{6}$$

In operator notation the above equations look like the following.

$$\left(1 - \Delta t \frac{\partial^2}{\partial x^2} - \Delta t \frac{\partial^2}{\partial y^2}\right) \psi^{n+1} = \psi^n + \Delta t \omega^n \tag{7}$$

$$\left(1 - \Delta t \frac{1}{\text{Re}} \frac{\partial^2}{\partial x^2} - \Delta t \frac{1}{\text{Re}} \frac{\partial^2}{\partial y^2} - \Delta t \left(\frac{\partial \psi}{\partial x}\right)^n \frac{\partial}{\partial y} + \Delta t \left(\frac{\partial \psi}{\partial y}\right)^n \frac{\partial}{\partial x}\right) \omega^{n+1} = \omega^n \tag{8}$$

The above equations are in fully implicit form where each equation requires the solution of a large banded matrix which is not computationally efficient. Instead of solving the equations (7 and 8) in fully implicit form, more efficiently, we spatially factorize the equations, such as

$$\left(1 - \Delta t \frac{\partial^2}{\partial x^2}\right)\left(1 - \Delta t \frac{\partial^2}{\partial y^2}\right) \psi^{n+1} = \psi^n + \Delta t \omega^n \tag{9}$$

$$\left(1 - \Delta t \frac{1}{\text{Re}} \frac{\partial^2}{\partial x^2} + \Delta t \left(\frac{\partial \psi}{\partial y}\right)^n \frac{\partial}{\partial x}\right)\left(1 - \Delta t \frac{1}{\text{Re}} \frac{\partial^2}{\partial y^2} - \Delta t \left(\frac{\partial \psi}{\partial x}\right)^n \frac{\partial}{\partial y}\right) \omega^{n+1} = \omega^n \tag{10}$$

Here, we split the single LHS operators into two operators, where each only contains derivatives in one direction. The only difference between equations (7 and 8) and (9 and 10) are the $\Delta t^2$ terms produced by the factorization. The advantage of this process is that each equation now require the solution of a tridiagonal system, which is numerically more efficient.

As the solution converges to the steady state solution, the $\Delta t^2$ terms due to the factorization will not become zero and will remain in the solution. To illustrate this, equation (9) can be written in explicit form as

$$\psi^{n+1} - \Delta t \frac{\partial^2 \psi^{n+1}}{\partial x^2} - \Delta t \frac{\partial^2 \psi^{n+1}}{\partial y^2} + \Delta t^2 \frac{\partial^4 \psi^{n+1}}{\partial x^2 \partial y^2} = \psi^n + \Delta t \omega^n \tag{11}$$

As the solution converges to a steady state, $\Psi^{n+1}$ becomes equal to $\Psi^n$ and they cancel each other from the above equation and at convergence, equation (11) therefore appears as the following

$$\frac{\partial^2 \psi}{\partial x^2} + \frac{\partial^2 \psi}{\partial y^2} + \omega = -\Delta t \frac{\partial^4 \psi}{\partial x^2 \partial y^2} \tag{12}$$

In equation (12), there is no guarantee that the RHS, which is the consequence of the factorization, will be small. Only a very small time step ($\Delta t$) will ensure that the RHS of equation (12) will be close to zero. This, however, will slow down the convergence. In addition, the product may not be small since $\frac{\partial^4 \psi}{\partial x^2 \partial y^2}$ can be $\mathcal{O}(\frac{1}{\Delta t})$ or higher in many flow situations.

In order to overcome this problem, we add a similar term at the previous time level to the RHS of equation (11) to give the following finite difference form

$$\psi^{n+1} - \Delta t \frac{\partial^2 \psi^{n+1}}{\partial x^2} - \Delta t \frac{\partial^2 \psi^{n+1}}{\partial y^2} + \Delta t^2 \frac{\partial^4 \psi^{n+1}}{\partial x^2 \partial y^2} = \psi^n + \Delta t \omega^n + \Delta t^2 \frac{\partial^4 \psi^n}{\partial x^2 \partial y^2} \tag{13}$$

As a result, at steady state, $\psi^{n+1}$ converges to $\psi^n$ and $\frac{\partial^4 \psi^{n+1}}{\partial x^2 \partial y^2}$ converges to $\frac{\partial^4 \psi^n}{\partial x^2 \partial y^2}$, so that they cancel each out, and the final solution converges to the correct physical form

$$\frac{\partial^2 \psi}{\partial x^2} + \frac{\partial^2 \psi}{\partial y^2} + \omega = 0 \tag{14}$$

We would like to emphasize that the added $\Delta t^2$ term on the RHS is not an artificial numerical diffusion term. Such terms are often used to stabilize the numerical schemes so that convergence can be achieved. Here, the added $\Delta t^2$ term on the RHS is not meant to stabilize the numerical scheme. The only reason for this approach is to cancel out the terms introduced by factorization so that the equations are the correct physical representation at the steady state.

One of the most popular spatially factorized schemes is undoubtedly the Beam & Warming [4] method. In that method, the equations are formulated in delta ($\Delta$) form. The main advantage of a delta formulation in a steady problem is that there will be no second order ($\Delta t^2$) terms due to the factorization in the solution at the steady state. However in our formulation the second order ($\Delta t^2$) terms due to factorization will remain even at the steady state. Adding a second order ($\Delta t^2$) term to the RHS of the equations to cancel out the terms produced by factorization, easily overcomes this problem. One of the main differences between the presented formulation and a delta formulation is in the numerical treatment of the vorticity equation. In a delta formulation of the vorticity equation, there appears $\Delta t$ order convection and dissipation terms on the RHS (explicit side) of the numerical equations. In the presented formulation, on the RHS of the equation there are no $\Delta t$ order terms. Instead cross derivative terms ($\frac{\partial^4}{\partial x^2 \partial y^2}$, $\frac{\partial^3}{\partial x^2 \partial y}$, $\frac{\partial^3}{\partial x \partial y^2}$, and $\frac{\partial^2}{\partial x \partial y}$) appear on the RHS. Note that these terms are $\Delta t^2$ order. Since $\Delta t$ is small, the contribution of these $\Delta t^2$ terms are even smaller. Our extensive numerical tests showed that this approach has better numerical stability characteristics and is more effective especially at high Reynolds numbers compared to a delta formulation of the steady streamfunction and vorticity equations.

Applying the approach of removing second order terms due to factorization on the RHS, to both equations (9 and 10), the final form of the numerical formulation we used for the solution of streamfunction and vorticity equation becomes

$$\left(1 - \Delta t \frac{\partial^2}{\partial x^2}\right)\left(1 - \Delta t \frac{\partial^2}{\partial y^2}\right) \psi^{n+1} = \psi^n + \Delta t \omega^n$$

$$+ \left(\Delta t \frac{\partial^2}{\partial x^2}\right)\left(\Delta t \frac{\partial^2}{\partial y^2}\right) \psi^n \tag{15}$$

$$\left(1 - \Delta t \frac{1}{\text{Re}} \frac{\partial^2}{\partial x^2} + \Delta t \left(\frac{\partial \psi}{\partial y}\right)^n \frac{\partial}{\partial x}\right)\left(1 - \Delta t \frac{1}{\text{Re}} \frac{\partial^2}{\partial y^2} - \Delta t \left(\frac{\partial \psi}{\partial x}\right)^n \frac{\partial}{\partial y}\right) \omega^{n+1} = \omega^n$$

$$+ \left(\Delta t \frac{1}{\text{Re}} \frac{\partial^2}{\partial x^2} - \Delta t \left(\frac{\partial \psi}{\partial y}\right)^n \frac{\partial}{\partial x}\right)\left(\Delta t \frac{1}{\text{Re}} \frac{\partial^2}{\partial y^2} + \Delta t \left(\frac{\partial \psi}{\partial x}\right)^n \frac{\partial}{\partial y}\right) \omega^n \tag{16}$$

The solution methodology for the two equations involves a two-level updating. For the streamfunction equation, the variable $f$ is introduced such that

$$\left(1 - \Delta t \frac{\partial^2}{\partial y^2}\right) \psi^{n+1} = f \tag{17}$$

and

$$\left(1 - \Delta t \frac{\partial^2}{\partial x^2}\right) f = \psi^n + \Delta t \omega^n + \left(\Delta t \frac{\partial^2}{\partial x^2}\right)\left(\Delta t \frac{\partial^2}{\partial y^2}\right) \psi^n \tag{18}$$

In equation (18) $f$ is the only unknown. First, this equation is solved for $f$ at each grid point. Following this, the streamfunction ($\psi$) is advanced into the new time level using equation (17).

In a similar fashion for the vorticity equation, the variable $g$ is introduced such that

$$\left(1 - \Delta t \frac{1}{\mathrm{Re}} \frac{\partial^2}{\partial y^2} - \Delta t \left(\frac{\partial \psi}{\partial x}\right)^n \frac{\partial}{\partial y}\right) \omega^{n+1} = g \tag{19}$$

and

$$\left(1 - \Delta t \frac{1}{\mathrm{Re}} \frac{\partial^2}{\partial x^2} + \Delta t \left(\frac{\partial \psi}{\partial y}\right)^n \frac{\partial}{\partial x}\right) g = \omega^n$$
$$+ \left(\Delta t \frac{1}{\mathrm{Re}} \frac{\partial^2}{\partial x^2} - \Delta t \left(\frac{\partial \psi}{\partial y}\right)^n \frac{\partial}{\partial x}\right)\left(\Delta t \frac{1}{\mathrm{Re}} \frac{\partial^2}{\partial y^2} + \Delta t \left(\frac{\partial \psi}{\partial x}\right)^n \frac{\partial}{\partial y}\right) \omega^n \tag{20}$$

As with $f$, $g$ is determined at every grid point using equation (20), then vorticity ($\omega$) is advanced into the next time level using equation (19).

We note that, in this numerical approach the streamfunction and vorticity equations are solved separately. Each equation is advanced into a next time level by solving two tridiagonal systems, which allows the use of very large grid sizes easily. The method proved to be very effective on flow problems that require high accuracy on very fine grid meshes (Erturk & Corke [9] Erturk, Haddad & Corke [10]).

## 3. RESULTS AND DISCUSSION

The boundary conditions and a schematics of the vortices generated in a driven cavity flow are shown in Figure 1. In this figure, the abbreviations BR, BL and TL refer to bottom right, bottom left and top left corners of the cavity, respectively. The number following these abbreviations refer to the vortices that appear in the flow, which are numbered according to size.

We have used the well-known Thom's formula [32] for the wall boundary condition, such that

$$\psi_0 = 0$$

$$\omega_0 = \frac{-2\psi_1}{\Delta h^2} - \frac{2U}{\Delta h} \tag{21}$$

where subscript 0 refers to points on the wall and 1 refers to points adjacent to the wall, $\Delta h$ refers to grid spacing and $U$ refers to the velocity of the wall with being equal to 1 on the moving wall and 0 on the stationary walls.

We note an important fact that, it is well understood ([34], [30], [18], [23]) that, even though Thom's method is locally first order accurate, the global solution obtained using Thom's method preserve second order accuracy. Therefore in this study, since three point second order central difference is used inside the cavity and Thom's method is used at the wall boundary conditions, the presented solutions are second order accurate $\mathcal{O}(\Delta h^2)$.

During our computations we monitored the residual of the steady streamfunction and vorticity equations (1 and 2) as a measure of the convergence to the steady state solution, where the residual of each equation is given as

$$R_\psi = \frac{\psi_{i-1,j}^{n+1} - 2\psi_{i,j}^{n+1} + \psi_{i+1,j}^{n+1}}{\Delta x^2} + \frac{\psi_{i,j-1}^{n+1} - 2\psi_{i,j}^{n+1} + \psi_{i,j+1}^{n+1}}{\Delta y^2} + \omega_{i,j}^{n+1} \qquad (22)$$

$$R_\omega = \frac{1}{\text{Re}} \frac{\omega_{i-1,j}^{n+1} - 2\omega_{i,j}^{n+1} + \omega_{i+1,j}^{n+1}}{\Delta x^2} + \frac{1}{\text{Re}} \frac{\omega_{i,j-1}^{n+1} - 2\omega_{i,j}^{n+1} + \omega_{i,j+1}^{n+1}}{\Delta y^2}$$

$$- \frac{\psi_{i,j+1}^{n+1} - \psi_{i,j-1}^{n+1}}{2\Delta y} \frac{\omega_{i+1,j}^{n+1} - \omega_{i-1,j}^{n+1}}{2\Delta x} + \frac{\psi_{i+1,j}^{n+1} - \psi_{i-1,j}^{n+1}}{2\Delta x} \frac{\omega_{i,j+1}^{n+1} - \omega_{i,j-1}^{n+1}}{2\Delta y} \qquad (23)$$

The magnitude of these residuals are an indication of the degree to which the solution has converged to steady state. In the limit these residuals would be zero. In our computations, for all Reynolds numbers, we considered that convergence was achieved when for each equation (22 and 23) the maximum of the absolute residual in the computational domain ($\max(|R_\psi|)$ and $\max(|R_\omega|)$) was less than $10^{-10}$. Such a low value was chosen to ensure the accuracy of the solution. At these residual levels, the maximum absolute difference in streamfunction value between two time steps, $(\max(|\psi^{n+1} - \psi^n|))$, was in the order of $10^{-16}$ and for vorticity, $(\max(|\omega^{n+1} - \omega^n|))$, it was in the order of $10^{-14}$. And also at these convergence levels, between two time steps the maximum absolute normalized difference in streamfunction, $(\max(|\frac{\psi^{n+1}-\psi^n}{\psi^n}|))$, and in vorticity, $(\max(|\frac{\omega^{n+1}-\omega^n}{\omega^n}|))$, were in the order of $10^{-13}$, and $10^{-12}$ respectively.

We have started our computations by solving the driven cavity flow from R=1,000 to Re=5,000 on a 129 × 129 grid mesh. With using this many number of grids, we could not get a steady solution for Re=7,500. The first natural conclusion was that the cavity flow may not be steady at Re=7,500, and therefore steady solution may not be computable. Even though we have used a pseudo-time iteration, the time evolution of the flow parameters, either vorticity or streamfunction values at certain locations, were periodic suggesting that the flow indeed may be periodic, backing up this conclusion. Besides, in the literature many studies claim that the driven cavity flow becomes unstable between a Reynolds number of 7,000 and 8,000, for example, Peng et. al. [25], Fortin et. al. [11] and Poliashenko & Aidun [26] predict unstable flow beyond Re of 7402, 7998,5 and 7763 respectively. However, on the other hand there were many studies ([2], [5], [12], [21], and [28]) that have presented solutions for a higher Reynolds number, Re=10,000. This was contradictory. We then have tried to solve the same case, Re=7,500, with a larger grid size of 257 × 257. This time we were able to obtain a steady solution. In fact, with 257 × 257 number of grids, we carry our steady computations up to Re=12,500. Barragy & Carey [2] have also presented steady solutions for Re=12,500. Interesting point was, when we tried to increase Re furthermore with the same grid mesh, we have again obtained a periodic behavior. Once again the natural conclusion was that the steady solutions beyond this Re may not be computable. Having concluded this, we could not help but ask the question what happens if we increase the number of grids furthermore. In their study Barragy & Carey [2] have presented a steady solution for Re=16,000. Their solution for this Reynolds number displayed oscillations. They concluded that this was due to a coarse mesh. This conclusion encouraged us to use larger grid size than $257 \times 257$ in order to obtain solutions for Reynolds numbers greater than 12,500. At this point, since the focus was on the number of grids, we realized that among the studies that presented steady solutions at very high Reynolds numbers ([2], [5], [12], [21], and [28]), the maximum number of grids used was $257 \times 257$. Larger grid sizes haven't been used at these high Reynolds numbers. This was also encouraging us to use more grid points. We, then, decided to increase the number of grids and tried to solve for Re=15,000 with $401 \times 401$ number of grids. This time we were able to obtain a steady solution. This fact suggests that in order to obtain a numerical solution for 2-D steady incompressible driven cavity flow for Re $\geq$ 12,500, a grid mesh with more than $257 \times 257$ grid points is needed. With using $401 \times 401$ grid points, we continued our steady computations up to Re=21,000. Beyond this Re our computations again displayed a periodic behavior. This time the first natural thing came to mind was to increase the number of grids. We have tried to use $513 \times 513$ grids, however again we could not get a steady solution beyond Re=21,000. Thinking that the increase in number of grids may not be enough, we then again increased the number of grids and have used $601 \times 601$ grids. The situation was the same, and the maximum Reynolds number that we can obtain a steady solution with using $601 \times 601$ grids was Re=21,000. We did not try to increase the number of grids furthermore since the computations became time consuming. Whether or not steady computations of driven cavity flow are possible beyond Re=21,000 with grid numbers larger than $601 \times 601$, is still an open question.

One of the reasons, why the steady solutions of the driven cavity flow at very high Reynolds numbers become computable when finer grids are used, may be the fact that as the number of grids used increases, $\Delta h$ gets smaller, then the cell Reynolds number or so called Peclet number defined as $\text{Re}_c = \frac{u \Delta h}{\nu}$ decreases. This improves the numerical stability characteristics of the numerical scheme used ([31], [34]), and allows high cavity Reynolds numbered solutions computable. Another reason may be that fact that finer grids would resolve the corner vortices better. This would, then, help decrease any numerical oscillations that might occur at the corners of the cavity during iterations.

We note that all the figures and tables present the solution of the finest grid size of $601 \times 601$ unless otherwise stated.

The accuracy of a finite difference solution is set by the mesh size, and by the spatial order of the finite difference equations and the boundary approximations. At low Reynolds numbers (Re=1,000), different numerical method solutions found in the literature agree with each other. However, the solutions at higher Reynolds numbers (Re>1,000) have noticeable discrepancies. We believe it is mainly due to different spatial orders and different grid mesh sizes and different boundary conditions used in different studies.

Table 1 tabulates the minimum streamfunction value, the vorticity value at the center of the primary vortex and also the center location of the primary vortex for Re=1,000 along with similar results found in the literature with the most significant ones are underlined. In Table 1 among the most significant (underlined) results, Schreiber & Keller [28] have used Richardson extrapolation in order to achieve high spatial accuracy. For Re=1,000 their extrapolated solutions are $\Delta h^6$ order accurate and these solutions are obtained by using repeated Richardson extrapolation on three different mesh size solutions.

Following the same approach, we have solved the cavity flow on three different grid mesh ($401 \times 401$, $513 \times 513$, and $601 \times 601$) separately for every Reynolds number (Re=1,000 to Re=21,000). Since the numerical solution is second order accurate $\mathcal{O}(\Delta h^2)$, the computed values of streamfunction and vorticity have an asymptotic error expansion of the form

$$\psi(x_i, y_j) = \psi_{i,j} + c_1 \Delta h^2 + c_2 \Delta h^4 + ...$$
$$\omega(x_i, y_j) = \omega_{i,j} + d_1 \Delta h^2 + d_2 \Delta h^4 + ... \tag{24}$$

We can use Richardson extrapolation in order to obtain high-order accurate approximations (Schreiber & Keller [28]). We have used repeated Richardson extrapolations using the three different mesh size solutions. The extrapolated values of the streamfunction and the vorticity at the center of the primary vortex are given in Table 2 and Table 3 respectively. The extrapolated results tabulated in Table 2 and Table 3 are, in theory, $\Delta h^6$ order accurate.

Looking back to Table 1, for Re=1,000, our extrapolated results are in very good agrement with the extrapolated results of Schreiber & Keller [28]. In fact our extrapolated streamfunction value at the primary vortex differs only 0.013 % and the vorticity value differs only 0.095 % from the extrapolated results of Schreiber & Keller [28].

With spectral methods, very high spatial accuracy can be obtained with a relatively smaller number of grid points. For this reason, the tabulated results from Botella & Peyret [6] are believed to be very accurate. Our results are also in very good agreement such that our extrapolated results differs 0.016 % in the streamfunction value and 0.098 % in the vorticity value from the results of Botella & Peyret [6].

Barragy & Carey [2] have only tabulated the location of the vortex centers and the streamfunction values at these locations yet they have not presented any vorticity data. Nevertheless, their solutions based on $p$-type finite element scheme also have very high spatial accuracy ($\Delta h^8$ order accurate). They have presented the maximum Reynolds number solutions (Re=12,500) for the cavity flow, to the best of our knowledge, found in the literature. We will show that the agreement between our results and results of Barragy & Carey [2] is excellent through the whole Reynolds number range they considered ($1,000 \leq \text{Re} \leq 12,500$). Comparing solutions for Re=1,000 in Table 1, the difference is 0.022 % in the streamfunction value.

Wright & Gaskell [36] have used a very fine grid mesh with $1024 \times 1024$ grid points in their solutions. With this many grids their solution should be quite accurate. Their streamfunction and vorticity value are very close to our solutions with a difference of 0.113 % and 0.114 % respectively.

Nishida & Satofuka [24] have presented higher order solutions in which for Re=1,000 their solutions are $\Delta h^8$ order accurate. Comparing with our results we find a very good agreement such that our extrapolated results differs 0.040 % in streamfunction value and 0.136 % in vorticity value.

Grigoriev & Dargush [14] have solved the governing equations with a boundary element method (BEM). They have only tabulated the streamfunction values and the location of the centers of the primary and secondary vortices for driven cavity flow for Reynolds numbers up to 5,000. Their solution for Re=1,000 is quite accurate.

Benjamin & Denny [5] have used a $\Delta h^n$ extrapolation to obtain the values when $\Delta h \rightarrow 0$ using the three different grid mesh size solutions. Their extrapolated results are close to other underlined data in Table 1. It is quite remarkable that even a simple extrapolation produces such good results.

Hou et. al. [17] have simulated the cavity flow by the Lattice Boltzmann Method. They have presented solutions for Re $\leq 7,500$ obtained by using $256 \times 256$ lattice. However their results have significant oscillations in the two upper corners of the cavity. Especially their vorticity contours have wiggles aligned with lattice directions, for Re $\geq$ 1,000.

Li et. al. [19] have used a compact scheme and presented fourth order accurate ($\Delta h^4$) results for a 129×129 grid mesh. Ghia et. al. [12] have used a second order method ($\Delta h^2$) with a mesh size of 129×129. These results, together with that of Goyon [13] and Liao & Zhu [21], might be considered as somewhat under-resolved. Our calculations showed that for a second order ($\Delta h^2$) spatial accuracy, even a 401 × 401 grid mesh solutions can be considered as under-resolved for Re=1,000.

The results of Bruneau & Jouron, [7] and Vanka [33] appear over-diffusive due to upwind differencing they have used.

From all these comparisons we can conclude that even for Re=1,000 higher order approximations together with the use of fine grids are necessary for accuracy.

Figures 2 through 6 show streamfunction contours of the finest (601 × 601) grid solutions, for various Reynolds numbers. These figures exhibit the formation of the counter-rotating secondary vortices which appear as the Reynolds number increases.

Table 4 compares our computed streamfunction and vorticity values at the center of the primary vortex with the results found in the literature for various Reynolds numbers. The agreement between our Richardson extrapolated results ($\Delta h^6$ order accurate) with results of Barragy & Carey [2] based on $p$-type finite element scheme ($\Delta h^8$ order accurate) is quite remarkable. One interesting result is that the streamfunction value at the center of the primary vortex decreases at first as the Reynolds number increases until Re=10,000, however then starts to increase as the Reynolds number increases further. This behavior was also documented by Barragy & Carey [2].

In terms of quantitative results, Gupta & Manohar [16], stated that for the numerical solutions of driven cavity flow, the value of the streamfunction ($\psi$) at the primary vortex center appears to be a reliable indicator of accuracy, where as the location of the vortex center ($x$ and $y$) is limited by the mesh size, and the value of the vorticity ($\omega$) at the vortex center is sensitive to the accuracy of the wall boundary conditions. Nevertheless, the locations of the primary and secondary vortices, as well as the values of the streamfunction ($\psi$) and vorticity ($\omega$) at these locations are presented in Table 5. These are comparable with the solutions of Barragy & Carey [2] and Ghia et. al. [12]. It is evident that as the Reynolds number increases, the center of the primary vortex moves towards the geometric center of the cavity. Its location is almost invariant for Re $\geq$ 17,500. Our computations indicate the appearance of a quaternary vortex at the bottom left corner (BL3) at Reynolds number of 10,000. We would like to note that, this quaternary vortex (BL3) did not appear in our solutions for Re=10,000 with grid sizes less than 513 × 513. This would suggest that, especially at high Reynolds numbers, in order to resolve the small vortices appear at the corners properly, very fine grids have to be used. Moreover, our solutions indicate another tertiary vortex at the top left side of the cavity (TL2) at Re=12,500. These vortices (BL3 and TL2) were also observed by Barragy & Carey [2].

Figures 7 and 8 present the u-velocity profiles along a vertical line and the v-velocity profiles along a horizontal line passing through the geometric center of the cavity respectively and also Tables 6 and 7 tabulates u- and v-velocity values at certain locations respectively, for future references. These profiles are in good agrement with that of Ghia et. al. [12] shown by symbols in Figures 7 and 8.

We have solved the incompressible flow in a driven cavity numerically. A very good mathematical check on the accuracy of the numerical solution would be to check the continuity of the fluid, as suggested by Aydin & Fenner [1], although this is done very rarely in driven cavity flow papers. We have integrated the u-velocity and v-velocity profiles, considered in Figures 7 and 8, along a vertical line and horizontal line passing through the geometric center of the cavity, in order to obtain the net volumetric flow rate, $Q$, through these sections. Since the flow is incompressible, the net volumetric flow rate passing through these sections should be equal to zero, $Q = 0$. The u-velocity and v-velocity profiles at every Reynolds number considered are integrated using Simpson's rule, then the obtained volumetric flow rate values are divided by a characteristic flow rate, $Q_c$, which is the horizontal rate that would occur in the absence of the side walls (Plane Couette flow), to help quantify the errors, as also suggested by Aydin & Fenner [1]. In an integration process the numerical errors will add up. Nevertheless, the volumetric flow rate values ($Q_1 = \frac{|\int_0^1 u dy|}{Q_c}$ and $Q_2 = \frac{|\int_0^1 v dx|}{Q_c}$) tabulated in Table 8 are close to zero such that, even the largest values of $Q_1 = 0.000000705$ and $Q_2 = 0.000002424$ in Table 8 can be considered as $Q_1 \approx Q_2 \approx 0$. This mathematical check on the conservation of the continuity shows that our numerical solution is indeed very accurate.

Distinctively than any other paper, Botella & Peyret [6] have tabulated highly accurate Chebyshev collocation spectral vorticity data from inside the cavity, along a vertical line and along a horizontal line passing through the geometric center of the cavity, for Re=1,000. Comparing our solutions with theirs ([6]) in Figure 9 and 10, we find that the agreement with Botella & Peyret [6] is remarkable, with the maximum difference between two solutions being 0.18 %.

Between the whole range of Re=1,000 and the maximum Re found in the literature, Re=12,500, our computed results agreed well with the published results. Since there is no previously published results (to the best of our knowledge) of steady driven cavity flow beyond Re > 12,500 in the literature, we decided to use Batchelor's [3] model in order to demonstrate the accuracy of our numerical solutions at higher Reynolds numbers. As seen in Figures 7 and 8, the u- and v-velocity profiles change almost linearly in the core of the primary vortex as Reynolds number increases. This would indicate that in this region the vorticity is uniform. As Re increases thin boundary layers are developed along the solid walls and the core fluid moves as a solid body with a uniform vorticity, in the manner suggested by Batchelor [3]. At high Reynolds numbers the vorticity at the core of the eddy is approximately constant and the flow in the core is governed by

$$\frac{\partial^2 \psi}{\partial x^2} + \frac{\partial^2 \psi}{\partial y^2} = C \qquad (25)$$

where $C$ is the constant vorticity value. With streamfunction value being $\psi = 0$ on the boundaries, inside the domain (x,y)=([-0.5,0.5],[-0.5,0.5]) the following expression is a solution to equation (25)

$$\psi = C \left( \frac{y^2}{2} - \frac{1}{8} + \sum_{n=1}^{\infty} A_n \cosh(2n-1)\pi x \cos(2n-1)\pi y \right) \qquad (26)$$

where

$$A_n = \frac{4(-1)^{n+1}}{(2n-1)^3 \pi^3} (\cosh(n - \frac{1}{2})\pi)^{-1} \qquad (27)$$

The mean square law proposed by Batchelor [3] relates the core vorticity value with the boundary velocities. From using the relation

$$\oint [(\psi_x)^2 + (\psi_y)^2] \mathrm{d}s = \oint |U|^2 \mathrm{d}s = 1 \tag{28}$$

where $U$ is the boundary velocity, we can have

$$\frac{1}{C^2} = \frac{1}{6} - \frac{16}{\pi^4} \sum_{n=1}^{\infty} \frac{1}{(2n-1)^4 (\cosh(n-\frac{1}{2})\pi)^2}$$

$$+ 2 \int_{-\frac{1}{2}}^{\frac{1}{2}} \left( \frac{1}{2} + \pi \sum_{n=1}^{\infty} (-1)^n (2n-1) A_n \cosh(2n-1)\pi x \right)^2 \mathrm{d}x \tag{29}$$

From the above expression, the numerical value of $C$ (vorticity value at the core) can be found approximately as 1.8859, in accordance with the theoretical core vorticity value of 1.886 at infinite Re, calculated analytically by Burggraf [8] in his application of Batchelor's model [3], which consists of an inviscid core with uniform vorticity, coupled to boundary layer flows at the solid surface.

In Figure 11 we have plotted the variation of the vorticity value at the center of the primary vortex, tabulated in Table 2, as a function of Reynolds number. In the same figure, the theoretical value of 1.8859 is shown with the dotted line. This figure clearly shows that our computations, in fact, in strong agrement with the mathematical theory (Batchelor [3], Wood [35], and Burggraf [8]) such that computed values of the vorticity at the primary vortex asymptote to the theoretical infinite Re vorticity value as the Reynolds number increases.

## 4. CONCLUSIONS

Fine grid numerical solutions of the steady driven cavity flow for Reynolds numbers up to Re=21,000 have been presented. The steady 2-D incompressible Navier-Stokes equations in streamfunction and vorticity formulation are solved computationally using the numerical method described. The streamfunction and vorticity equations are solved separately. For each equation, the numerical formulation requires the solution of two tridiagonal systems, which allows the use of large grid meshes easily. With this we have used a fine grid mesh of $601 \times 601$. At all Reynolds numbers, the numerical solutions converged to maximum absolute residuals of the governing equations that were less than $10^{-10}$.

Our computations indicate that fine grid mesh is necessary in order to obtain a steady solution and also resolve the vortices appear at the corners of the cavity, as the Reynolds number increases. Several unique features of the driven cavity flow have been presented such that the appearance of a quaternary vortex at the bottom left corner (BL3) at Reynolds number of 10,000 on fine grid mesh, and a tertiary vortex at the top left corner (TL2) at Reynolds number of 12,500. Detailed results were presented and they were compared with results found in the literature. The presented results were found to agree well with the published numerical solutions and with the theory as well.

ACKNOWLEDGEMENT

This study was funded by Gebze Institute of Technology with project no BAP-2003-A-22. E. Ertürk and C. Gökçöl are grateful for this financial support.

REFERENCES

1. M. Aydin, R.T. Fenner. Boundary Element Analysis of Driven Cavity Flow for Low and Moderate Reynolds Numbers. *Int. J. Numer. Meth. Fluids* 2001; **37**:45–64.
2. E. Barragy, G.F. Carey. Stream Function-Vorticity Driven Cavity Solutions Using *p* Finite Elements. *Computers and Fluids* 1997; **26**:453–468.
3. G.K. Batchelor. On Steady Laminar Flow with Closed Streamlines at Large Reynolds Numbers. *J. Fluid Mech.* 1956; **1**:177–190.
4. R.M. Beam, R.F. Warming. An Implicit Factored Scheme for the Compressible Navier-Stokes Equations. *AIAA J.* 1978; **16**:393–420.
5. A.S. Benjamin, V.E. Denny. On the Convergence of Numerical Solutions for 2-D Flows in a Cavity at Large Re. *J. Comp. Physics* 1979; **33**:340–358.
6. O. Botella, R. Peyret. Benchmark Spectral Results on the Lid-Driven Cavity Flow. *Computers and Fluids* 1998; **27**:421–433.
7. C-H. Bruneau, C. Jouron. An Efficient Scheme for Solving Steady Incompressible Navier-Stokes Equations. *J. Comp. Physics* 1990; **89**:389–413.
8. O.R. Burggraf. Analytical and Numerical Studies of the Structure of Steady Separated Flows. *J. Fluid Mech.* 1966; **24**:113–151.
9. E. Erturk, T.C. Corke. Boundary layer leading-edge receptivity to sound at incidence angles. *J. Fluid Mech.* 2001; **444**:383–407.
10. E. Erturk, O.M. Haddad, T.C. Corke. Laminar Incompressible Flow Past Parabolic Bodies at Angles of Attack. *AIAA J.* 2004; **42**:2254–2265.
11. A. Fortin, M. Jardak, J.J. Gervais, R. Pierre. Localization of Hopf Bifurcations in Fluid Flow Problems. *Int. J. Numer. Methods Fluids* 1997; **24**:1185–1210.
12. U. Ghia, K.N. Ghia, C.T. Shin. High-Re Solutions for Incompressible Flow Using the Navier-Stokes Equations and a Multigrid Method. *J. Comp. Physics* 1982; **48**:387–411.
13. O. Goyon. High-Reynolds Number Solutions of Navier-Stokes Equations Using Incremental Unknowns. *Computer Methods in Applied Mechanics and Engineering* 1996; **130**:319–335.
14. M.M. Grigoriev, G.F. Dargush. A Poly-Region Boundary Element Method for Incompressible Viscous Fluid Flows. *Int. J. Numer. Meth. Engng.* 1999; **46**:1127–1158.
15. M.M. Gupta. High Accuracy Solutions of Incompressible Navier-Stokes Equations. *J. Comp. Physics* 1991; **93**:343–359.
16. M.M. Gupta, R.P. Manohar. Boundary Approximations and Accuracy in Viscous Flow Computations. *J. Comp. Physics* 1979; **31**:265–288.
17. S. Hou, Q. Zou, S. Chen, G. Doolen, A.C. Cogley. Simulation of Cavity Flow by the Lattice Boltzmann Method. *J. Comp. Physics* 1995; **118**:329–347.
18. H. Huang, B.R. Wetton. Discrete Compatibility in Finite Difference Methods for Viscous Incompressible Fluid Flow. *J. Comp. Physics* 1996; **126**:468–478.
19. M. Li, T. Tang, B. Fornberg. A Compact Forth-Order Finite Difference Scheme for the Steady Incompressible Navier-Stokes Equations. *Int. J. Numer. Methods Fluids* 1995; **20**:1137–1151.
20. S.J. Liao. Higher-Order Streamfunction-Vorticity Formulation of 2-D Steady State Navier-Stokes Equations. *Int. J. Numer. Methods Fluids* 1992; **15**:595–612.
21. S.J. Liao, J.M. Zhu. A Short Note on Higher-Order Stremfunction-Vorticity Formulation of 2-D Steady State Navier-Stokes Equations. *Int. J. Numer. Methods Fluids* 1996; **22**:1–9.
22. M. Nallasamy, K.K. Prasad. On Cavity Flow at High Reynolds Numbers. *J. Fluid Mech.* 1977; **79**:391–414.
23. M. Napolitano, G. Pascazio, L. Quartapelle. A Review of Vorticity Conditions in the Numerical Solution of the $\zeta$-$\psi$ Equations. *Computers and Fluids* 1999; **28**:139–185.
24. H. Nishida, N. Satofuka. Higher-Order Solutions of Square Driven Cavity Flow Using a Variable-Order Multi-Grid Method. *Int. J. Numer. Methods Fluids* 1992; **34**:637–653.
25. Y.F. Peng, Y.H. Shiau, R.R. Hwang, Transition in a 2-D Lid-Driven Cavity Flow. *Computers & Fluids* 2003; **32**:337–352.
26. M. Poliashenko, C.K. Aidun. A Direct Method for Computation of Simple Bifurcations. *Journal of Computational Physics* 1995; **121**:246–260.

27. S.G. Rubin, P.K. Khosla. Navier-Stokes Calculations with a Coupled Strongly Implicit Method. *Computers and Fluids* 1981; **9**:163–180.
28. R. Schreiber, H.B. Keller. Driven Cavity Flows by Efficient Numerical Techniques. *J. Comp. Physics* 1983; **49**:310–333.
29. P.N. Shankar, M.D. Deshpande. Fluid Mechanics in the Driven Cavity. *Annu. Rev. Fluid Mech.* 2000; **32**:93–136.
30. W.F. Spotz. Accuracy and Performance of Numerical Wall Boundary Conditions for Steady 2D Incompressible Streamfunction Vorticity. *Int. J. Numer. Methods Fluids* 1998; **28**:737–757.
31. J.C. Tennehill, D.A. Anderson, R.H. Pletcher. *Computational Fluid Mechanics and Heat Transfer* (2nd edn). Taylor & Francis, 1997.
32. A. Thom. The Flow Past Circular Cylinders at Low Speed. *Proc. Roy. Soc. Lond. A.* 1933; **141**:651–669
33. S.P. Vanka. Block-Implicit Multigrid Solution of Navier-Stokes Equations in Primitive Variables. *J. Comp. Physics* 1986; **65**:138–158.
34. E. Weinan, L. Jian-Guo. Vorticity Boundary Condition and Related Issues for Finite Difference Schemes. *J. Comp. Physics* 1996; **124**:368–382.
35. W.W. Wood. Boundary Layers Whose Streamlines are Closed. *J. Fluid Mech.* 1957; **2**:77–87.
36. N.G. Wright, P.H. Gaskell. An efficient Multigrid Approach to Solving Highly Recirculating Flows. *Computers and Fluids* 1995; **24**:63–79.

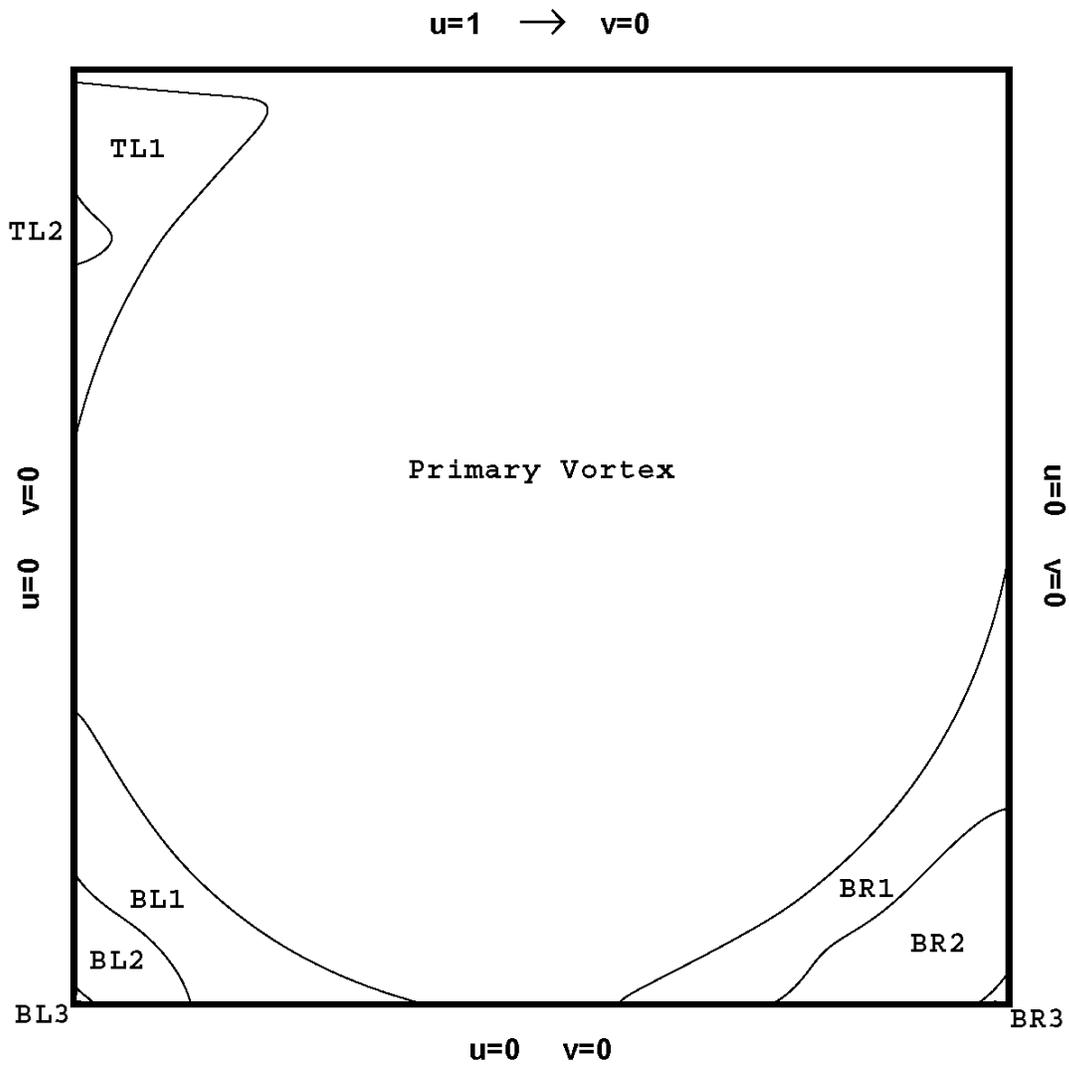

Figure 1. Schematic view of driven cavity flow

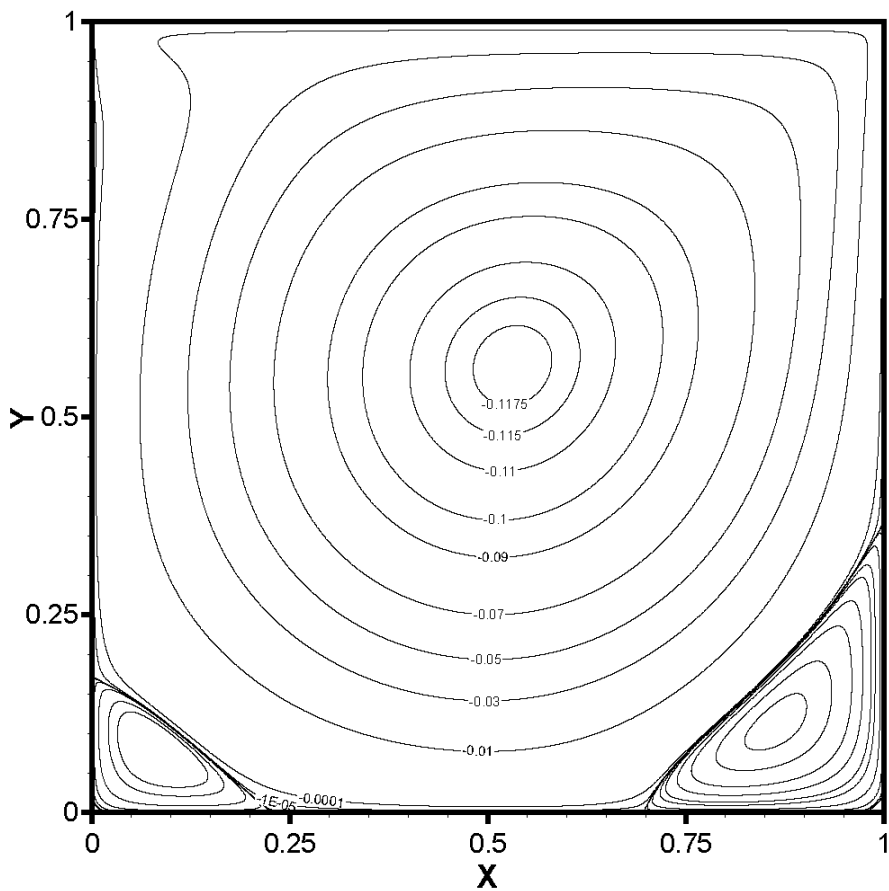

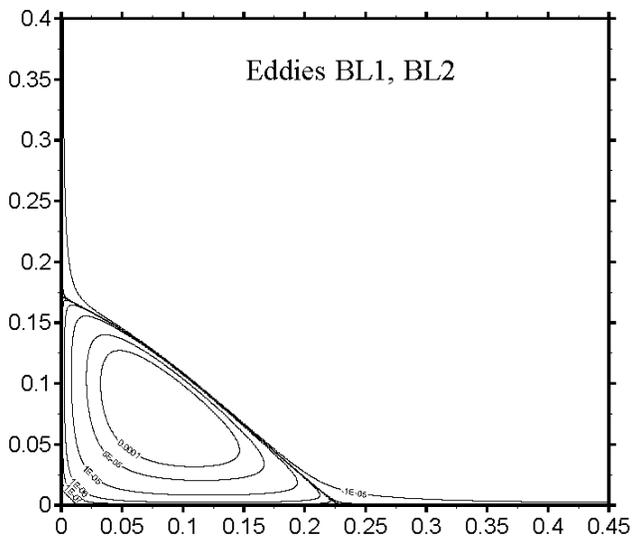 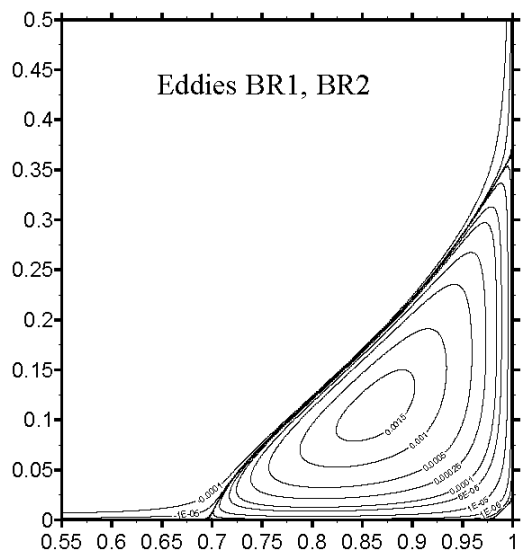

Figure 2. Streamline contours of primary and secondary vortices, Re=1,000

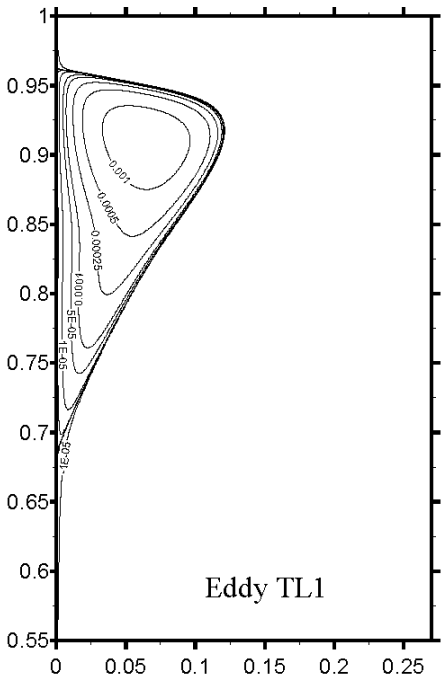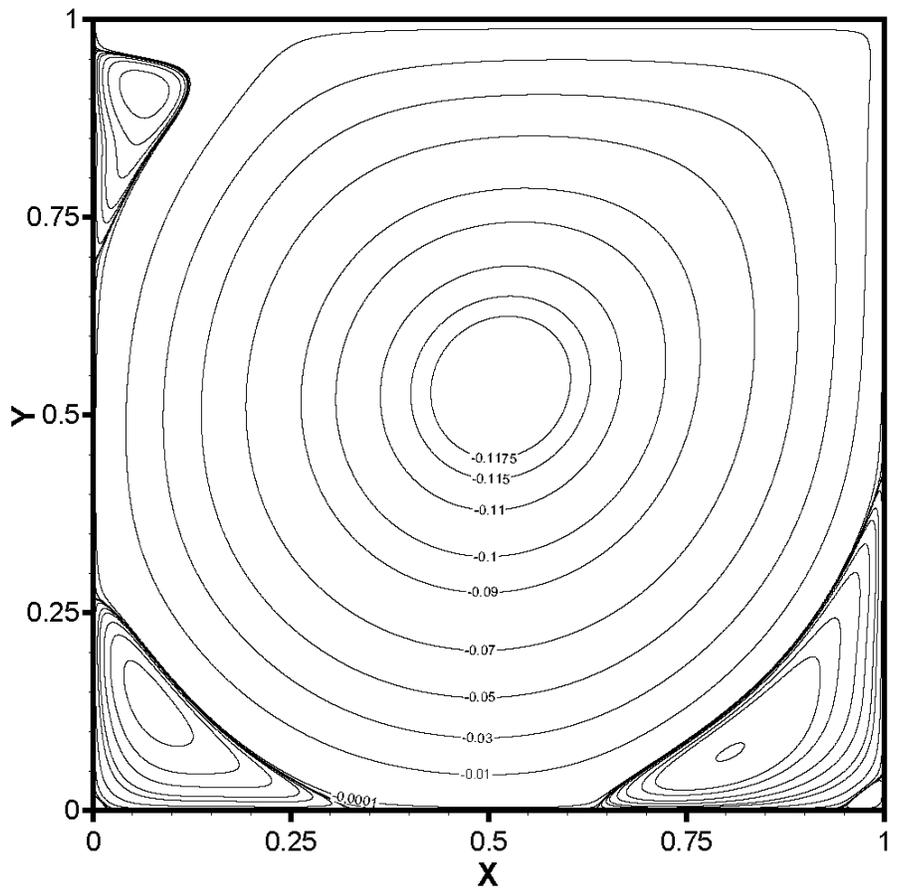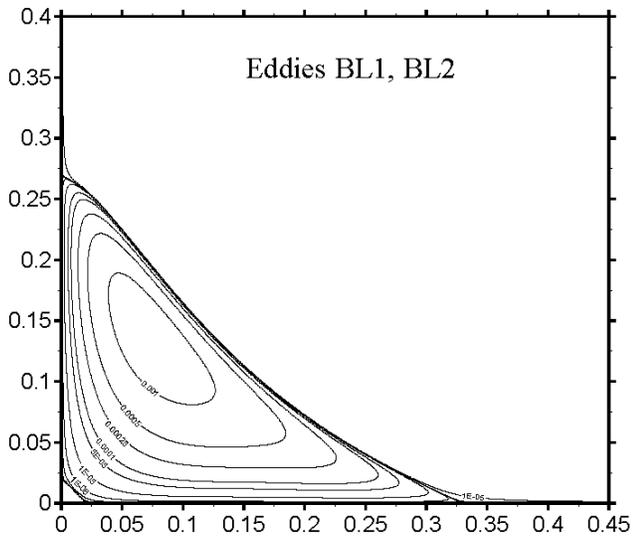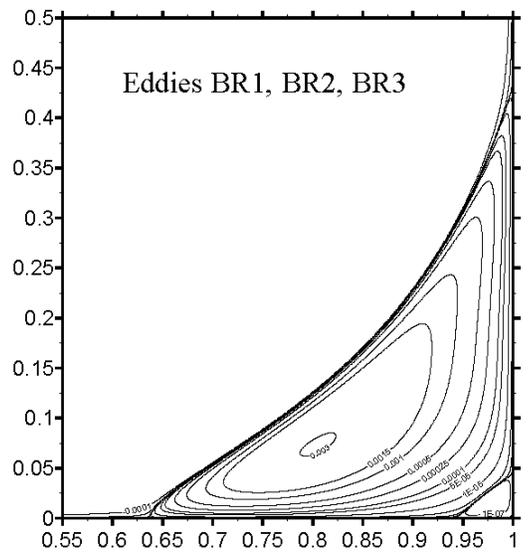

Figure 3. Streamline contours of primary and secondary vortices, Re=5,000

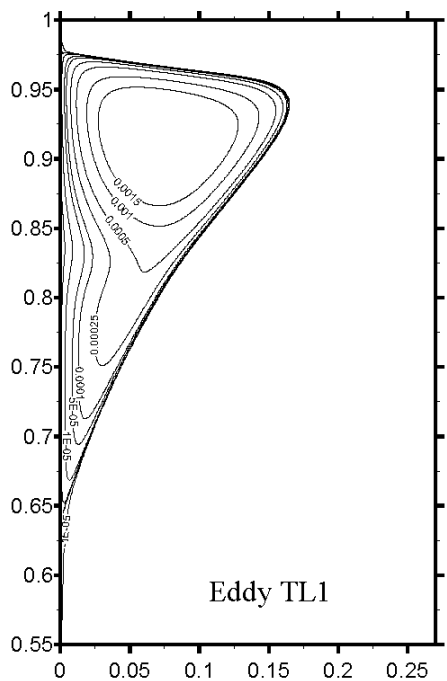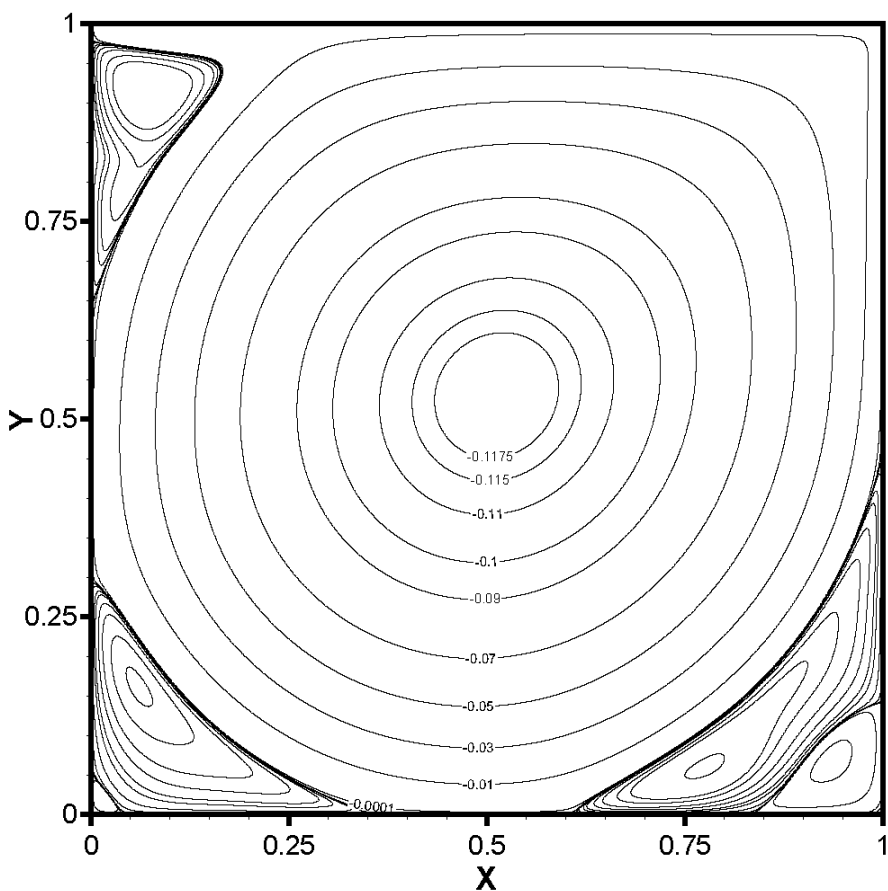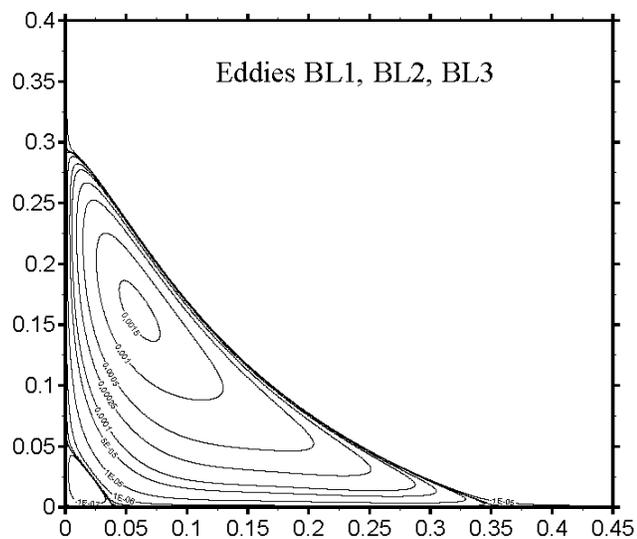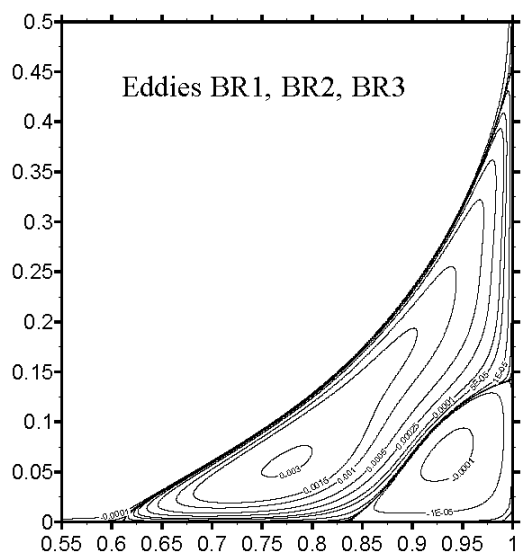

Figure 4. Streamline contours of primary and secondary vortices, Re=10,000

Figure 5. Streamline contours of primary and secondary vortices, Re=15,000

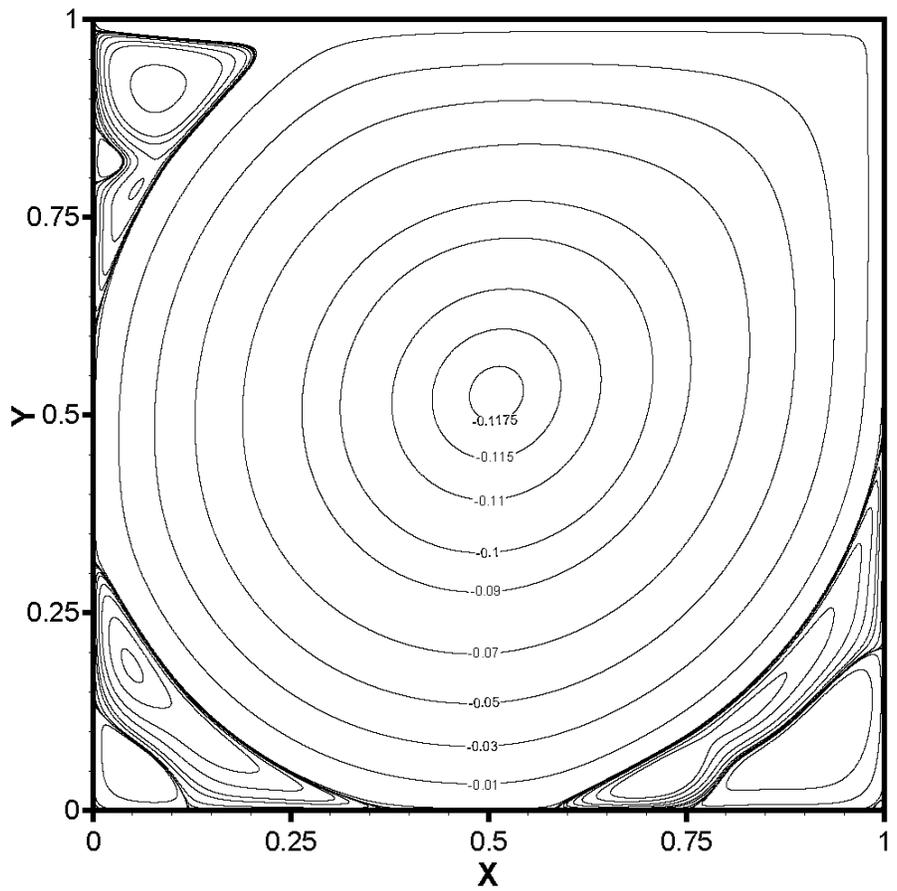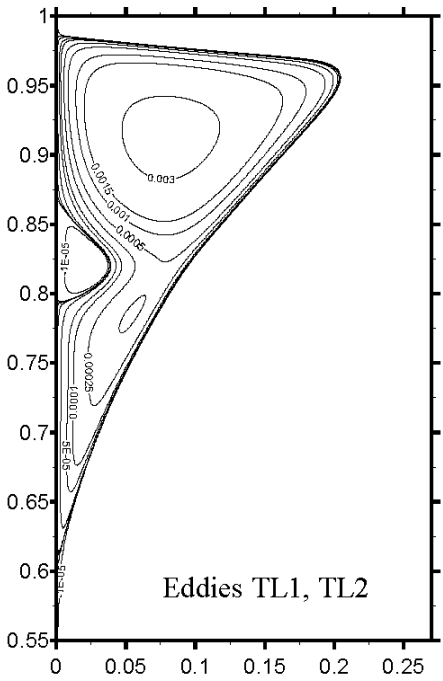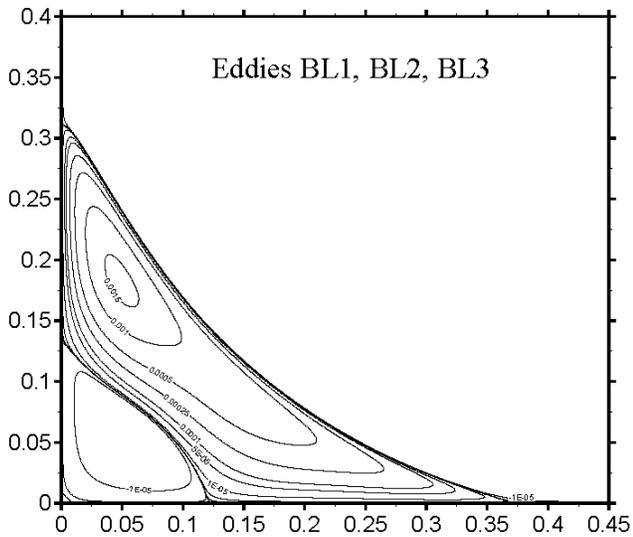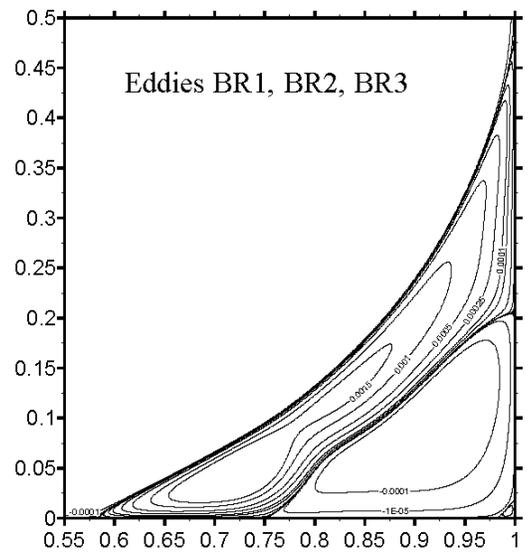

Figure 6. Streamline contours of primary and secondary vortices, Re=20,000

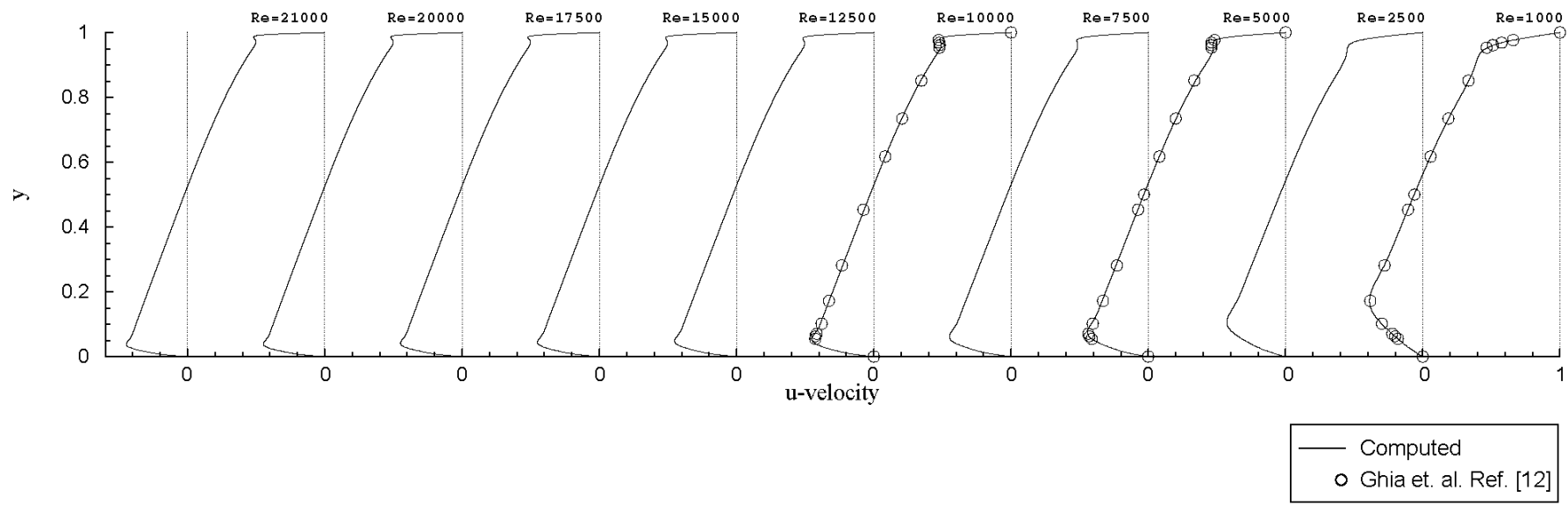

Figure 7. Computed u-velocity profiles profiles along a vertical line passing through the geometric center of the cavity at various Reynolds numbers

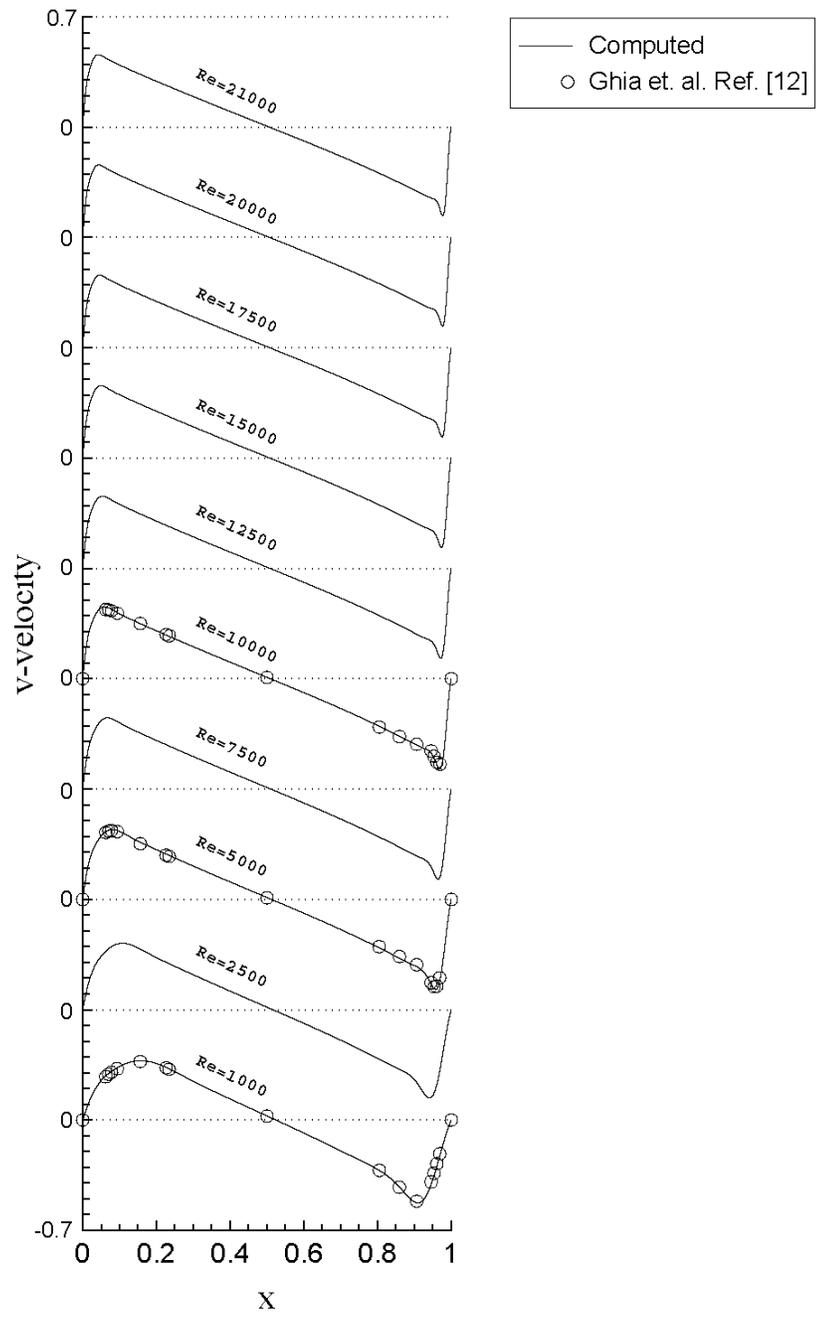

Figure 8. Computed v-velocity profiles along a horizontal line passing through the geometric center of the cavity at various Reynolds numbers

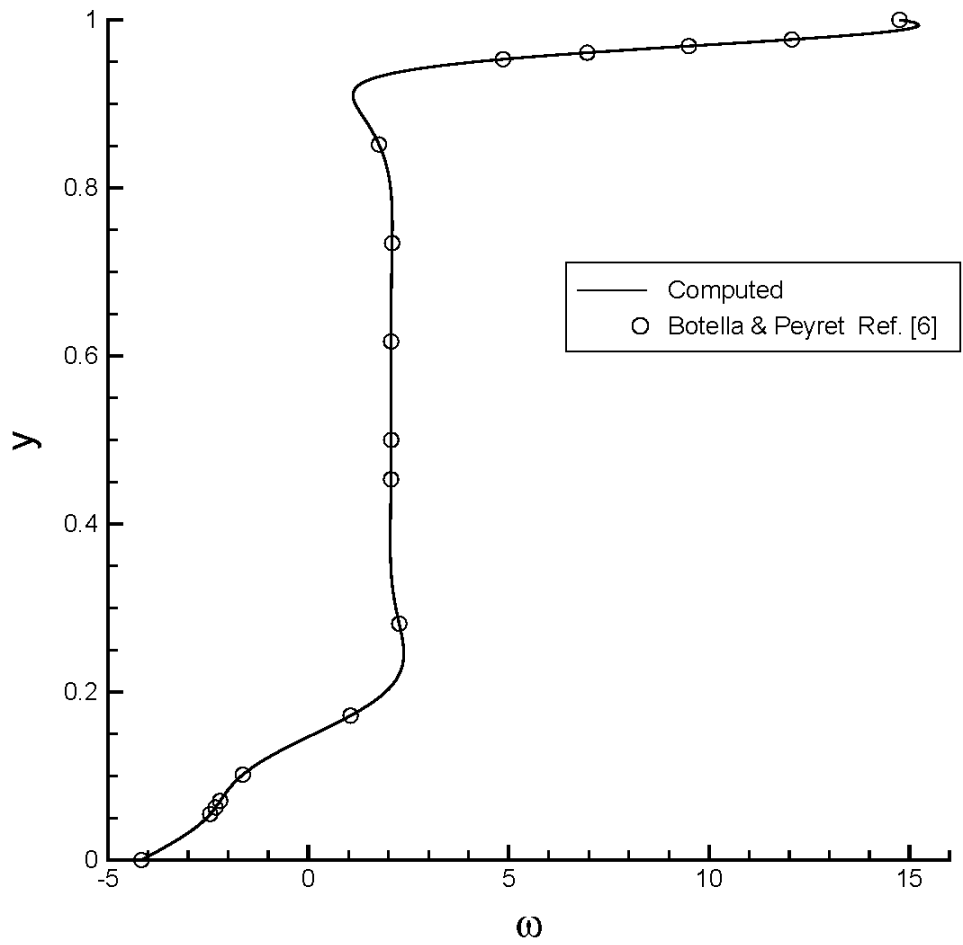

Figure 9.  Computed vorticity values along a vertical line passing through the geometric center of the cavity, Re=1,000

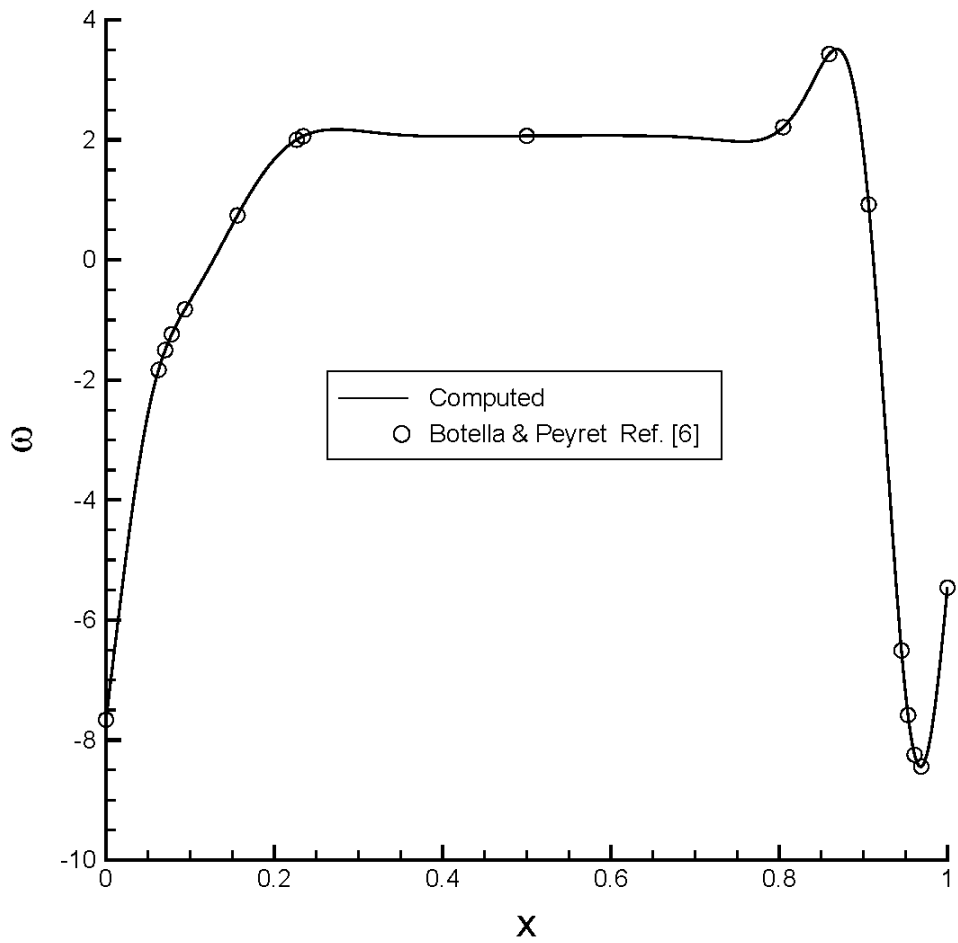

Figure 10. Computed vorticity values along a horizontal line passing through the geometric center of the cavity, Re=1,000

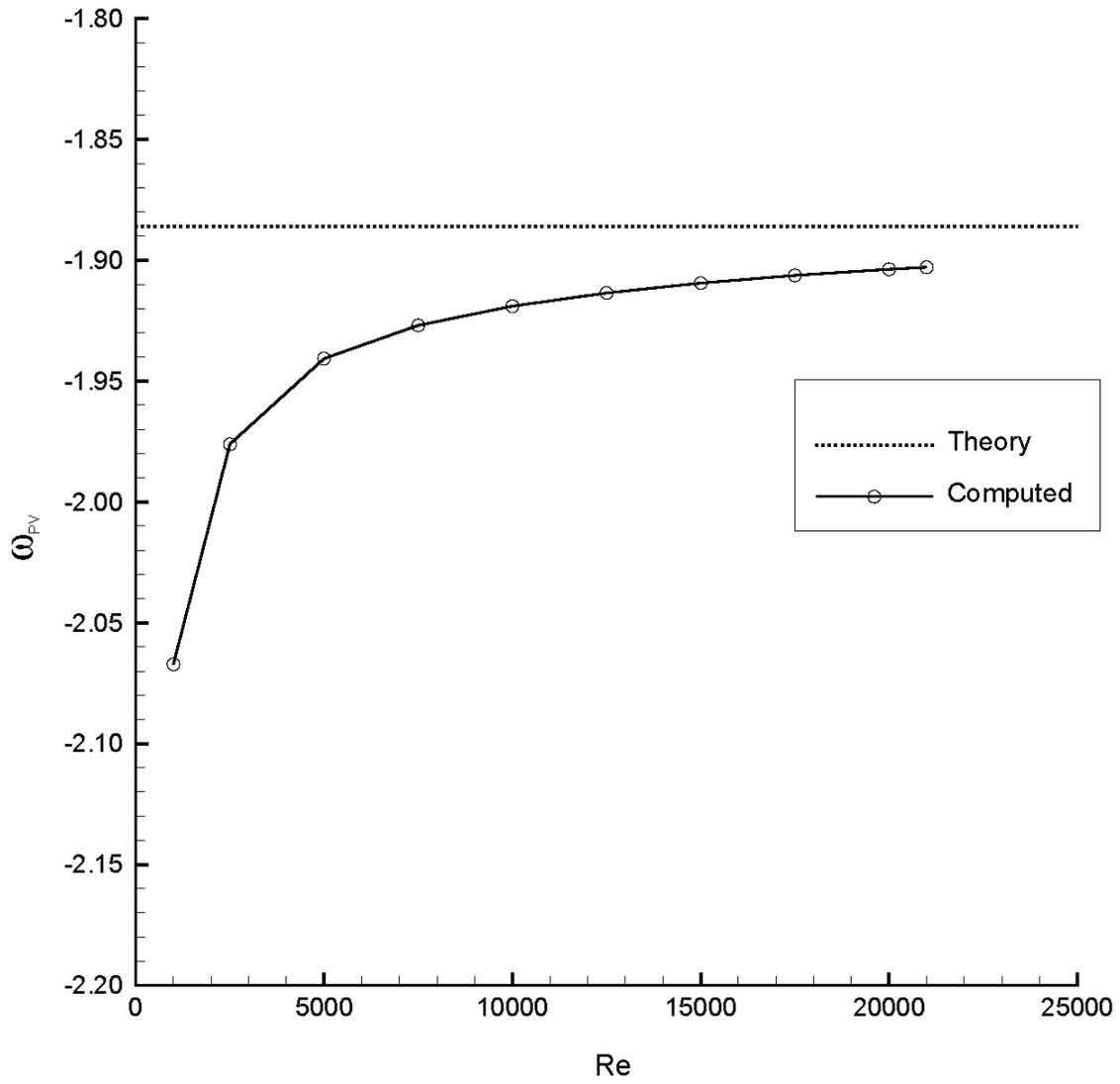

Figure 11. Strength of vorticity at the primary vortex as a function of Reynolds number.

| Reference | No | Grid | Spatial Accuracy | $\psi$ | $\omega$ | x | y |
|---|---|---|---|---|---|---|---|
| Present | | 401 × 401 | $\Delta h^2$ | 0.118585 | 2.062761 | 0.5300 | 0.5650 |
| Present | | 513 × 513 | $\Delta h^2$ | 0.118722 | 2.064765 | 0.5313 | 0.5645 |
| Present | | 601 × 601 | $\Delta h^2$ | 0.118781 | 2.065530 | 0.5300 | 0.5650 |
| Present | | Extrapolated | $\Delta h^6$ | 0.118942 | 2.067213 | - | - |
| Barragy & Carey | Ref. [2] | 257 × 257 | $p$=8 | 0.118930 | - | - | - |
| Botella & Peyret | Ref. [6] | N=128 | N=128 | 0.1189366 | 2.067750 | 0.5308 | 0.5652 |
| Botella & Peyret | Ref. [6] | N=160 | N=160 | 0.1189366 | 2.067753 | 0.5308 | 0.5652 |
| Schreiber & Keller | Ref. [26] | 100 × 100 | $\Delta h^2$ | 0.11315 | 1.9863 | - | - |
| Schreiber & Keller | Ref. [26] | 121 × 121 | $\Delta h^2$ | 0.11492 | 2.0112 | - | - |
| Schreiber & Keller | Ref. [26] | 141 × 141 | $\Delta h^2$ | 0.11603 | 2.0268 | 0.52857 | 0.56429 |
| Schreiber & Keller | Ref. [26] | Extrapolated | $\Delta h^6$ | 0.11894 | 2.0677 | - | - |
| Wright & Gaskell | Ref. [34] | 1024 × 1024 | $\Delta h^2$ | 0.118821 | 2.06337 | 0.5308 | 0.5659 |
| Nishida & Satofuka | Ref. [22] | 129 × 129 | $\Delta h^8$ | 0.119004 | 2.068546 | 0.5313 | 0.5625 |
| Benjamin & Denny | Ref. [5] | 101 × 101 | $\Delta h^2$ | 0.1175 | 2.044 | - | - |
| Benjamin & Denny | Ref. [5] | Extrapolated | high order | 0.1193 | 2.078 | - | - |
| Li et. al. | Ref. [18] | 129 × 129 | $\Delta h^4$ | 0.118448 | 2.05876 | 0.5313 | 0.5625 |
| Ghia et. al. | Ref. [12] | 129 × 129 | $\Delta h^2$ | 0.117929 | 2.04968 | 0.5313 | 0.5625 |
| Bruneau & Jouron | Ref. [7] | 256 × 256 | $\Delta h^2$ | 0.1163 | - | 0.5313 | 0.5586 |
| Goyon | Ref. [13] | 129 × 129 | $\Delta h^2$ | 0.1157 | - | 0.5312 | 0.5625 |
| Vanka | Ref. [31] | 321 × 321 | $\Delta h^2$ | 0.1173 | - | 0.5438 | 0.5625 |
| Gupta | Ref. [15] | 41 × 41 | $\Delta h^4$ | 0.111492 | 2.02763 | 0.525 | 0.575 |
| Hou et. al. | Ref. [17] | 256 × 256 | $\Delta h^2$ | 0.1178 | 2.0760 | 0.5333 | 0.5647 |
| Liao & Zhu | Ref. [20] | 129 × 129 | $\Delta h^2$ | 0.1160 | 2.0234 | 0.5313 | 0.5625 |
| Grigoriev & Dargush | Ref. [14] | - | - | 0.11925 | - | 0.531 | 0.566 |

Table 1. Comparison of the properties of the primary vortex; the maximum streamfunction value, the vorticity value and the location of the center, for Re=1000

| Re | Grid | $\psi\ \mathcal{O}(\Delta h^2)$ | $\psi\ \mathcal{O}(\Delta h^4)$ | $\psi\ \mathcal{O}(\Delta h^6)$ |
|---|---|---|---|---|
| 1000 | 401×401 | -0.118585 | | |
| | 513×513 | -0.118722 | -0.118937 | |
| | 601×601 | -0.118781 | -0.118939 | -0.118942 |
| 2500 | 401×401 | -0.120493 | | |
| | 513×513 | -0.120873 | -0.121468 | |
| | 601×601 | -0.121035 | -0.121469 | -0.121470 |
| 5000 | 401×401 | -0.120145 | | |
| | 513×513 | -0.120944 | -0.122196 | |
| | 601×601 | -0.121289 | -0.122213 | -0.122233 |
| 7500 | 401×401 | -0.119178 | | |
| | 513×513 | -0.120395 | -0.122301 | |
| | 601×601 | -0.120924 | -0.122341 | -0.122386 |
| 10000 | 401×401 | -0.118059 | | |
| | 513×513 | -0.119690 | -0.122245 | |
| | 601×601 | -0.120403 | -0.122313 | -0.122390 |
| 12500 | 401×401 | -0.116889 | | |
| | 513×513 | -0.118936 | -0.122142 | |
| | 601×601 | -0.119831 | -0.122229 | -0.122326 |
| 15000 | 401×401 | -0.115701 | | |
| | 513×513 | -0.118162 | -0.122017 | |
| | 601×601 | -0.119239 | -0.122124 | -0.122245 |
| 17500 | 401×401 | -0.114507 | | |
| | 513×513 | -0.117381 | -0.121883 | |
| | 601×601 | -0.118641 | -0.122016 | -0.122167 |
| 20000 | 401×401 | -0.113313 | | |
| | 513×513 | -0.116596 | -0.121739 | |
| | 601×601 | -0.118038 | -0.121901 | -0.122084 |
| 21000 | 401×401 | -0.112837 | | |
| | 513×513 | -0.116282 | -0.121678 | |
| | 601×601 | -0.117797 | -0.121855 | -0.122056 |

Table 2. Richardson extrapolation for streamfunction values at the primary vortex

| Re | Grid | $\omega\ \mathcal{O}(\Delta h^2)$ | $\omega\ \mathcal{O}(\Delta h^4)$ | $\omega\ \mathcal{O}(\Delta h^6)$ |
|---|---|---|---|---|
| 1000 | 401×401 | -2.062761 | | |
| | 513×513 | -2.064765 | -2.067904 | |
| | 601×601 | -2.065530 | -2.067579 | -2.067213 |
| 2500 | 401×401 | -1.961660 | | |
| | 513×513 | -1.967278 | -1.976078 | |
| | 601×601 | -1.969675 | -1.976096 | -1.976117 |
| 5000 | 401×401 | -1.909448 | | |
| | 513×513 | -1.921431 | -1.940201 | |
| | 601×601 | -1.926601 | -1.940451 | -1.940732 |
| 7500 | 401×401 | -1.878523 | | |
| | 513×513 | -1.896895 | -1.925673 | |
| | 601×601 | -1.904883 | -1.926282 | -1.926969 |
| 10000 | 401×401 | -1.853444 | | |
| | 513×513 | -1.878187 | -1.916945 | |
| | 601×601 | -1.888987 | -1.917919 | -1.919018 |
| 12500 | 401×401 | -1.830890 | | |
| | 513×513 | -1.862010 | -1.910757 | |
| | 601×601 | -1.875618 | -1.912072 | -1.913557 |
| 15000 | 401×401 | -1.809697 | | |
| | 513×513 | -1.847199 | -1.905943 | |
| | 601×601 | -1.863618 | -1.907602 | -1.909476 |
| 17500 | 401×401 | -1.789348 | | |
| | 513×513 | -1.833212 | -1.901921 | |
| | 601×601 | -1.852447 | -1.903975 | -1.906294 |
| 20000 | 401×401 | -1.769594 | | |
| | 513×513 | -1.819761 | -1.898343 | |
| | 601×601 | -1.841814 | -1.900891 | -1.903767 |
| 21000 | 401×401 | -1.761828 | | |
| | 513×513 | -1.814492 | -1.896986 | |
| | 601×601 | -1.837672 | -1.899768 | -1.902909 |

Table 3. Richardson extrapolation for vorticity values at the primary vortex

|  |  |  | Re |  |  |  |  |
|--|--|--|--|--|--|--|--|
|  |  |  | 2500 | 5000 | 7500 | 10000 | 12500 |
| Present |  | $\psi$ | -0.121470 | -0.122233 | -0.122386 | -0.122390 | -0.122326 |
|  |  | $\omega$ | -1.976117 | -1.940732 | -1.926969 | -1.919018 | -1.913557 |
| Barragy & Carey | Ref. [2] | $\psi$ | -0.1214621 | -0.122219 | -0.1223803 | -0.1223930 | -0.1223584 |
|  |  | $\omega$ | - | - | - | - | - |
| Schreiber & Keller | Ref. [26] | $\psi$ | - | - | - | -0.12292 | - |
|  |  | $\omega$ | - | - | - | -1.9263 | - |
| Li et. al. | Ref. [18] | $\psi$ | - | -0.120359 | -0.1193791 | - | - |
|  |  | $\omega$ | - | -1.92430 | -1.91950 | - | - |
| Benjamin & Denny | Ref. [5] | $\psi$ | - | - | - | -0.1212 | - |
|  |  | $\omega$ | - | - | - | - | - |
| Ghia et. al. | Ref. [12] | $\psi$ | - | -0.118966 | -0.119976 | -0.119731 | - |
|  |  | $\omega$ | - | -1.86016 | -1.87987 | -1.88082 | - |
| Liao & Zhu | Ref. [20] | $\psi$ | - | -0.1186 | -0.1201 | -0.1201 | - |
|  |  | $\omega$ | - | -1.9375 | -1.9630 | -1.9426 | - |
| Hou et. al. | Ref. [17] | $\psi$ | - | -0.1214 | -0.1217 | - | - |
|  |  | $\omega$ | - | - | - | - | - |
| Grigoriev & Dargush | Ref. [14] | $\psi$ | - | -0.12209 | - | - | - |
|  |  | $\omega$ | - | - | - | - | - |

Table 4. Comparison of the minimum streamfunction value, the vorticity value at the center of the primary vortex for various Reynolds numbers

| | Re | 1000 | 2500 | 5000 | 7500 | 10000 | 12500 | 15000 | 17500 | 20000 | 21000 |
|---|---|---|---|---|---|---|---|---|---|---|---|
| **Primary** | $\psi$ | -0.118781 | -0.121035 | -0.121289 | -0.120924 | -0.120403 | -0.119831 | -0.119240 | -0.118641 | -0.118039 | -0.117797 |
| **Vortex** | $\omega$ | -2.065530 | -1.969675 | -1.926601 | -1.904883 | -1.888987 | -1.875618 | -1.863618 | -1.852447 | -1.841814 | -1.837672 |
| | (x , y) | (0.5300 , 0.5650) | (0.5200 , 0.5433) | (0.5150 , 0.5350) | (0.5133 , 0.5317) | (0.5117 , 0.5300) | (0.5117 , 0.5283) | (0.5100 , 0.5283) | (0.5100 , 0.5267) | (0.5100 , 0.5267) | (0.5100 , 0.5267) |
| | $\psi$ | 0.17281E-02 | 0.26561E-02 | 0.30604E-02 | 0.32102E-02 | 0.31758E-02 | 0.30839E-02 | 0.29897E-02 | 0.28938E-02 | 0.27967E-02 | 0.27598E-02 |
| **BR1** | $\omega$ | 1.115505 | 1.929029 | 2.724481 | 3.215956 | 3.696991 | 4.260493 | 4.837680 | 5.372747 | 5.882504 | 6.080586 |
| | (x , y) | (0.8633 , 0.1117) | (0.8350 , 0.0917) | (0.8050 , 0.0733) | (0.7900 , 0.0650) | (0.7767 , 0.0600) | (0.7617 , 0.0550) | (0.7483 , 0.0500) | (0.7367 , 0.0467) | (0.7267 , 0.0450) | (0.7217 , 0.0433) |
| | $\psi$ | 0.23261E-03 | 0.92541E-03 | 0.13639E-02 | 0.15171E-02 | 0.15930E-02 | 0.16355E-02 | 0.16484E-02 | 0.16331E-02 | 0.16032E-02 | 0.15897E-02 |
| **BL1** | $\omega$ | 0.353473 | 0.966506 | 1.502628 | 1.838156 | 2.150847 | 2.318733 | 2.441871 | 2.594288 | 2.823121 | 2.877385 |
| | (x , y) | (0.0833 , 0.0783) | (0.0850 , 0.1100) | (0.0733 , 0.1367) | (0.0650 , 0.1517) | (0.0583 , 0.1633) | (0.0550 , 0.1683) | (0.0533 , 0.1717) | (0.0517 , 0.1750) | (0.0483 , 0.1817) | (0.0483 , 0.1833) |
| | $\psi$ | -0.54962E-07 | -0.12967E-06 | -0.14010E-05 | -0.30998E-04 | -0.13397E-03 | -0.24687E-03 | -0.33148E-03 | -0.39539E-03 | -0.44959E-03 | -0.47028E-03 |
| **BR2** | $\omega$ | -0.77076E-02 | -0.15959E-01 | -0.34108E-01 | -0.158392 | -0.302615 | -0.397507 | -0.453339 | -0.512637 | -0.544236 | -0.574871 |
| | (x , y) | (0.9917 , 0.0067) | (0.9900 , 0.0100) | (0.9783 , 0.0183) | (0.9517 , 0.0417) | (0.9350 , 0.0667) | (0.9283 , 0.0817) | (0.9267 , 0.0883) | (0.9283 , 0.0967) | (0.9300 , 0.1033) | (0.9317 , 0.1083) |
| | $\psi$ | -0.84221E-08 | -0.32321E-07 | -0.72080E-07 | -0.20202E-06 | -0.10151E-05 | -0.57409E-05 | -0.20172E-04 | -0.44887E-04 | -0.74543E-04 | -0.86415E-04 |
| **BL2** | $\omega$ | -0.29802E-02 | -0.93800E-02 | -0.12252E-01 | -0.16862E-01 | -0.30981E-01 | -0.68557E-01 | -0.134120 | -0.192907 | -0.232067 | -0.250751 |
| | (x , y) | (0.0050 , 0.0050) | (0.0067 , 0.0067) | (0.0083 , 0.0083) | (0.0117 , 0.0117) | (0.0167 , 0.0200) | (0.0250 , 0.0317) | (0.0367 , 0.0417) | (0.0483 , 0.0483) | (0.0567 , 0.0533) | (0.0600 , 0.0550) |
| | $\psi$ | - | 0.34455E-03 | 0.14416E-02 | 0.21119E-02 | 0.25870E-02 | 0.29394E-02 | 0.32110E-02 | 0.34298E-02 | 0.36105E-02 | 0.36749E-02 |
| **TL1** | $\omega$ | - | 1.321269 | 2.062775 | 2.201598 | 2.247599 | 2.296811 | 2.329198 | 2.358627 | 2.380582 | 2.387858 |
| | (x , y) | - | (0.0433 , 0.8900) | (0.0633 , 0.9100) | (0.0667 , 0.9133) | (0.0717 , 0.9117) | (0.0750 , 0.9117) | (0.0767 , 0.9117) | (0.0800 , 0.9133) | (0.0817 , 0.9133) | (0.0817 , 0.9133) |
| | $\psi$ | - | - | 0.21562E-09 | 0.14409E-08 | 0.51934E-08 | 0.99032E-08 | 0.13947E-07 | 0.19492E-07 | 0.29610E-07 | 0.34929E-07 |
| **BR3** | $\omega$ | - | - | 0.70828E-03 | 0.22180E-02 | 0.25744E-02 | 0.38075E-02 | 0.33588E-02 | 0.43303E-02 | 0.51574E-02 | 0.63746E-02 |
| | (x , y) | - | - | (0.9983 , 0.0017) | (0.9967 , 0.0033) | (0.9967 , 0.0050) | (0.9950 , 0.0050) | (0.9950 , 0.0050) | (0.9950 , 0.0067) | (0.9933 , 0.0067) | (0.9933 , 0.0083) |
| | $\psi$ | - | - | - | - | 0.14382E-09 | 0.50930E-09 | 0.11280E-08 | 0.22898E-08 | 0.34584E-08 | 0.37846E-08 |
| **BL3** | $\omega$ | - | - | - | - | 0.78229E-03 | 0.53920E-03 | 0.10824E-02 | 0.22499E-02 | 0.18414E-02 | 0.17128E-02 |
| | (x , y) | - | - | - | - | (0.0017 , 0.0017) | (0.0017 , 0.0017) | (0.0017 , 0.0033) | (0.0033 , 0.0033) | (0.0033 , 0.0033) | (0.0033 , 0.0033) |
| | $\psi$ | - | - | - | - | - | -0.12752E-05 | -0.13394E-04 | -0.36218E-04 | -0.63084E-04 | -0.74273E-04 |
| **TL2** | $\omega$ | - | - | - | - | - | -0.215534 | -0.515658 | -0.711567 | -0.858460 | -0.933479 |
| | (x , y) | - | - | - | - | - | (0.0067 , 0.8317) | (0.0150 , 0.8267) | (0.0200 , 0.8233) | (0.0233 , 0.8200) | (0.0250 , 0.8200) |

Table 5. Properties of primary and secondary vortices; streamfunction and vorticity values, (x,y) locations

|       | Re     |        |        |        |        |        |        |        |        |        |
| y     | 1000   | 2500   | 5000   | 7500   | 10000  | 12500  | 15000  | 17500  | 20000  | 21000  |
| ---   | ---    | ---    | ---    | ---    | ---    | ---    | ---    | ---    | ---    | ---    |
| 1.000 | 1.0000 | 1.0000 | 1.0000 | 1.0000 | 1.0000 | 1.0000 | 1.0000 | 1.0000 | 1.0000 | 1.0000 |
| 0.990 | 0.8486 | 0.7704 | 0.6866 | 0.6300 | 0.5891 | 0.5587 | 0.5358 | 0.5183 | 0.5048 | 0.5003 |
| 0.980 | 0.7065 | 0.5924 | 0.5159 | 0.4907 | 0.4837 | 0.4833 | 0.4850 | 0.4871 | 0.4889 | 0.4895 |
| 0.970 | 0.5917 | 0.4971 | 0.4749 | 0.4817 | 0.4891 | 0.4941 | 0.4969 | 0.4982 | 0.4985 | 0.4983 |
| 0.960 | 0.5102 | 0.4607 | 0.4739 | 0.4860 | 0.4917 | 0.4937 | 0.4937 | 0.4925 | 0.4906 | 0.4897 |
| 0.950 | 0.4582 | 0.4506 | 0.4738 | 0.4824 | 0.4843 | 0.4833 | 0.4811 | 0.4784 | 0.4754 | 0.4742 |
| 0.940 | 0.4276 | 0.4470 | 0.4683 | 0.4723 | 0.4711 | 0.4684 | 0.4653 | 0.4622 | 0.4592 | 0.4580 |
| 0.930 | 0.4101 | 0.4424 | 0.4582 | 0.4585 | 0.4556 | 0.4523 | 0.4492 | 0.4463 | 0.4436 | 0.4425 |
| 0.920 | 0.3993 | 0.4353 | 0.4452 | 0.4431 | 0.4398 | 0.4366 | 0.4338 | 0.4312 | 0.4287 | 0.4277 |
| 0.910 | 0.3913 | 0.4256 | 0.4307 | 0.4275 | 0.4243 | 0.4216 | 0.4190 | 0.4166 | 0.4142 | 0.4132 |
| 0.900 | 0.3838 | 0.4141 | 0.4155 | 0.4123 | 0.4095 | 0.4070 | 0.4047 | 0.4024 | 0.4001 | 0.3992 |
| 0.500 | -0.0620 | -0.0403 | -0.0319 | -0.0287 | -0.0268 | -0.0256 | -0.0247 | -0.0240 | -0.0234 | -0.0232 |
| 0.200 | -0.3756 | -0.3228 | -0.3100 | -0.3038 | -0.2998 | -0.2967 | -0.2942 | -0.2920 | -0.2899 | -0.2892 |
| 0.180 | -0.3869 | -0.3439 | -0.3285 | -0.3222 | -0.3179 | -0.3146 | -0.3119 | -0.3096 | -0.3074 | -0.3066 |
| 0.160 | -0.3854 | -0.3688 | -0.3467 | -0.3406 | -0.3361 | -0.3326 | -0.3297 | -0.3271 | -0.3248 | -0.3239 |
| 0.140 | -0.3690 | -0.3965 | -0.3652 | -0.3587 | -0.3543 | -0.3506 | -0.3474 | -0.3446 | -0.3422 | -0.3412 |
| 0.120 | -0.3381 | -0.4200 | -0.3876 | -0.3766 | -0.3721 | -0.3685 | -0.3652 | -0.3622 | -0.3595 | -0.3585 |
| 0.100 | -0.2960 | -0.4250 | -0.4168 | -0.3978 | -0.3899 | -0.3859 | -0.3827 | -0.3797 | -0.3769 | -0.3758 |
| 0.080 | -0.2472 | -0.3979 | -0.4419 | -0.4284 | -0.4142 | -0.4054 | -0.4001 | -0.3965 | -0.3936 | -0.3925 |
| 0.060 | -0.1951 | -0.3372 | -0.4272 | -0.4491 | -0.4469 | -0.4380 | -0.4286 | -0.4206 | -0.4143 | -0.4121 |
| 0.040 | -0.1392 | -0.2547 | -0.3480 | -0.3980 | -0.4259 | -0.4407 | -0.4474 | -0.4490 | -0.4475 | -0.4463 |
| 0.020 | -0.0757 | -0.1517 | -0.2223 | -0.2633 | -0.2907 | -0.3113 | -0.3278 | -0.3412 | -0.3523 | -0.3562 |
| 0.000 | 0.0000 | 0.0000 | 0.0000 | 0.0000 | 0.0000 | 0.0000 | 0.0000 | 0.0000 | 0.0000 | 0.0000 |

Table 6. Tabulated u-velocity profiles along a vertical line passing through the geometric center of the cavity at various Reynolds numbers

| x | Re 1000 | 2500 | 5000 | 7500 | 10000 | 12500 | 15000 | 17500 | 20000 | 21000 |
|---|---|---|---|---|---|---|---|---|---|---|
| 1.000 | 0.0000 | 0.0000 | 0.0000 | 0.0000 | 0.0000 | 0.0000 | 0.0000 | 0.0000 | 0.0000 | 0.0000 |
| 0.985 | -0.0973 | -0.1675 | -0.2441 | -0.2991 | -0.3419 | -0.3762 | -0.4041 | -0.4269 | -0.4457 | -0.4522 |
| 0.970 | -0.2173 | -0.3725 | -0.5019 | -0.5550 | -0.5712 | -0.5694 | -0.5593 | -0.5460 | -0.5321 | -0.5266 |
| 0.955 | -0.3400 | -0.5192 | -0.5700 | -0.5434 | -0.5124 | -0.4899 | -0.4754 | -0.4664 | -0.4605 | -0.4588 |
| 0.940 | -0.4417 | -0.5603 | -0.5139 | -0.4748 | -0.4592 | -0.4534 | -0.4505 | -0.4482 | -0.4459 | -0.4449 |
| 0.925 | -0.5052 | -0.5268 | -0.4595 | -0.4443 | -0.4411 | -0.4388 | -0.4361 | -0.4331 | -0.4300 | -0.4287 |
| 0.910 | -0.5263 | -0.4741 | -0.4318 | -0.4283 | -0.4256 | -0.4221 | -0.4186 | -0.4153 | -0.4122 | -0.4110 |
| 0.895 | -0.5132 | -0.4321 | -0.4147 | -0.4118 | -0.4078 | -0.4040 | -0.4005 | -0.3975 | -0.3946 | -0.3936 |
| 0.880 | -0.4803 | -0.4042 | -0.3982 | -0.3938 | -0.3895 | -0.3859 | -0.3828 | -0.3800 | -0.3774 | -0.3764 |
| 0.865 | -0.4407 | -0.3843 | -0.3806 | -0.3755 | -0.3715 | -0.3682 | -0.3654 | -0.3627 | -0.3603 | -0.3593 |
| 0.850 | -0.4028 | -0.3671 | -0.3624 | -0.3574 | -0.3538 | -0.3508 | -0.3481 | -0.3457 | -0.3434 | -0.3425 |
| 0.500 | 0.0258 | 0.0160 | 0.0117 | 0.0099 | 0.0088 | 0.0080 | 0.0074 | 0.0069 | 0.0065 | 0.0063 |
| 0.150 | 0.3756 | 0.3918 | 0.3699 | 0.3616 | 0.3562 | 0.3519 | 0.3483 | 0.3452 | 0.3423 | 0.3413 |
| 0.135 | 0.3705 | 0.4078 | 0.3878 | 0.3779 | 0.3722 | 0.3678 | 0.3641 | 0.3608 | 0.3579 | 0.3567 |
| 0.120 | 0.3605 | 0.4187 | 0.4070 | 0.3950 | 0.3885 | 0.3840 | 0.3801 | 0.3767 | 0.3736 | 0.3725 |
| 0.105 | 0.3460 | 0.4217 | 0.4260 | 0.4137 | 0.4056 | 0.4004 | 0.3964 | 0.3929 | 0.3897 | 0.3885 |
| 0.090 | 0.3273 | 0.4142 | 0.4403 | 0.4337 | 0.4247 | 0.4180 | 0.4132 | 0.4093 | 0.4060 | 0.4048 |
| 0.075 | 0.3041 | 0.3950 | 0.4426 | 0.4495 | 0.4449 | 0.4383 | 0.4323 | 0.4273 | 0.4232 | 0.4218 |
| 0.060 | 0.2746 | 0.3649 | 0.4258 | 0.4494 | 0.4566 | 0.4563 | 0.4529 | 0.4484 | 0.4438 | 0.4420 |
| 0.045 | 0.2349 | 0.3238 | 0.3868 | 0.4210 | 0.4409 | 0.4522 | 0.4580 | 0.4602 | 0.4601 | 0.4596 |
| 0.030 | 0.1792 | 0.2633 | 0.3263 | 0.3608 | 0.3844 | 0.4018 | 0.4152 | 0.4254 | 0.4332 | 0.4357 |
| 0.015 | 0.1019 | 0.1607 | 0.2160 | 0.2509 | 0.2756 | 0.2940 | 0.3083 | 0.3197 | 0.3290 | 0.3323 |
| 0.000 | 0.0000 | 0.0000 | 0.0000 | 0.0000 | 0.0000 | 0.0000 | 0.0000 | 0.0000 | 0.0000 | 0.0000 |

Table 7. Tabulated v-velocity profiles along a horizontal line passing through the geometric center of the cavity at various Reynolds numbers

| Re | $Q_1 = \frac{|\int_0^1 u\,dy|}{Q_c}$ | $Q_2 = \frac{|\int_0^1 v\,dx|}{Q_c}$ |
|---|---|---|
| 1000 | 0.000000045 | 0.000000134 |
| 2500 | 0.000000046 | 0.000000344 |
| 5000 | 0.000000067 | 0.000000693 |
| 7500 | 0.000000090 | 0.000001027 |
| 10000 | 0.000000114 | 0.000001341 |
| 12500 | 0.000000143 | 0.000001630 |
| 15000 | 0.000000163 | 0.000001894 |
| 17500 | 0.000000259 | 0.000002133 |
| 20000 | 0.000000557 | 0.000002346 |
| 21000 | 0.000000705 | 0.000002424 |

Table 8. Volumetric flow rates through a vertical line, $Q_1$, and a horizontal line, $Q_2$, passing through the geometric center of the cavity